\documentclass[10pt,twoside,a4paper,reqno]{amsart}
\usepackage{amsmath}
\usepackage{amssymb}
\renewcommand{\k}{\kappa}
\newcommand{\sm}{\setminus}
\newcommand{\om}{\omega}
\newcommand{\R}{\mathbb{R}}
\newcommand{\su}{\subseteq}

\newcommand{\rest}{\restriction}

\newcommand{\frow}{{}^\frown }

\newcommand{\bd}{\operatorname{bd}}
\newcommand{\supt}{\operatorname{supt}}

\newcommand{\Lip}{\mathcal Lip}
\newcommand{\cont}{\mathcal Cont}
\newcommand{\Cont}{\mathcal Cont}  %%%% Alert: I added \cont=\Cont

\newcommand{\cl}{\operatorname{cl}}

\newcommand{\cf}{\operatorname{cf}}

\newcommand{\norm}{\operatorname{norm}}

\newcommand{\dom}{\operatorname{dom}}
\newcommand{\comp}{\operatorname{comp}}
\newcommand{\Comp}{\operatorname{Comp}}

\renewcommand{\subsetneq}{\subsetneqq}

\newcommand{\cov}{\operatorname{Cov}}

\renewcommand{\int}{\operatorname{int}}

\newcommand{\dist}{\operatorname{dist}}

\newcommand{\id}{\operatorname{id}}

\renewcommand{\succ}{\operatorname{succ}}
\newcommand{\Succ}{\operatorname{succ}}
\newcommand{\hm}{\mathfrak{hm}}

\newcommand{\parity}{\operatorname{parity}}

\newcommand{\random}{\operatorname{random}}

\renewcommand{\d}{\mathfrak{d}}

\newcommand{\forces}{\Vdash}

\newcommand{\concatenate}{\frow}

\newcommand{\card}[1]{\mathopen{|}{#1}\mathclose{|}}

 \newcommand{\hsqueeze}[1]{\hbox to 0pt{\hss$#1$\hss}}
  \newcommand{\nop}{}   % [1]{\hbox to 0pt{\hss$#1$\hss}}

\newcommand{\LineNotArrow}{
  \tmptmp=\LineLength
  \advance\tmptmp\RoomAboveLine
  \vrule width 0.1pt height \LineLength depth 0cm
  \vrule width 0pt   height \tmptmp     depth  \RoomAboveLine }

\newdimen\LineLength
\newdimen\RoomAboveLine
\newdimen\tmptmp

\newtheorem*{ab}{Theorem}
\newtheorem{thm}{Theorem}[section]
\newtheorem{corollary}[thm]{Corollary}
\newtheorem{lemma}[thm]{Lemma}
\newtheorem{problem}[thm]{Problem}

\newtheorem{defn}[thm]{Definition}
\newtheorem{claim}[thm]{Claim}
\newtheorem{rem}[thm]{Remark}

\newtheorem{fact}[thm]{Fact}
\newtheorem{Notation}[thm]{Notation}

\begin{document}

\title[Continuous Ramsey theory]{Continuous  Ramsey theory on
Polish spaces and covering the plane by functions}
\author{Stefan Geschke}
\author{Martin Goldstern} \author{Menachem Kojman}

\address[Geschke]{II.~Mathematisches Institut\\Freie Universit\"at
Berlin\\Arnimallee 3\\14195 Berlin\\Germany}
\email{geschke@math.fu-berlin.de}
\address[Goldstern]{Algebra\\TU Wien\\
Wiedner Hauptstrasse 8-10/118\\A-1040 Wien\\
Austria, Europe}
\email{goldstern@tuwien.ac.at}
\address[Kojman]{Department of Mathematics\\Ben Gurion University
of the Negev\\Beer Sheva, Israel}
\email{kojman@math.bgu.ac.il}
\date{May 15, 2002}
\subjclass[2000]{Primary: 03E17, 03E35, 05C55; Secondary: 26A15, 26A16}
\keywords{pair coloring, continuous coloring, homogeneity number,
Polish space, Cantor set, Lipschitz function, covering number, optimal forcing,
tree forcing}

\begin{abstract} We investigate the Ramsey theory of continuous pair-colorings
on complete, separable metric spaces, and apply the results to
the problem of covering a plane by functions.

The \emph{homogeneity number $\hm(c)$} of a pair-coloring $c:[X]^2\to
2$ is the number of $c$-homogeneous subsets of $X$ needed to
cover~$X$. We isolate two continuous pair-colorings on the Cantor
space $2^\om$, $c_{\min}$ and $c_{\max}$, which satisfy
$\hm(c_{\min})\le \hm(c_{\max})$ and prove:

\begin{ab}
\begin{enumerate}
\item For every
Polish space
$X$ and    every continuous pair-coloring $c:[X]^2\to 2$ with
$\hm(c)>\aleph_0$,

\[\hm(c)=\hm(c_{\min}) \text {\; \emph{ or }\; } \hm(c)=\hm(c_{\max}).\]

\item There is  a model of set theory in which
$\hm(c_{\min})=\aleph_1$ and $\hm(c_{\max})=\aleph_2$.
\end{enumerate}
\end{ab}

The
consistency of $\hm(c_{\min})=2^{\aleph_0}$ and of
$\hm(c_{\max})<2^{\aleph_0}$ follows from \cite{GKKS}.

We  prove that $\hm(c_{\min})$ is equal to the covering number of
$(2^\om)^2$ by graphs of Lipschitz functions and their
reflections on the diagonal. An iteration of an optimal forcing
notion associated to $c_{\min}$ gives:

\begin{ab} There is a model of set theory in
which
\begin{enumerate}
\item $\R^2$ is coverable by $\aleph_1$ graphs and reflections of
graphs of \emph{continuous} real functions;
\item $\R^2$ is not coverable by
$\aleph_1$ graphs and reflections of graphs  of \emph{Lipschitz} real
functions.
\end{enumerate}
\end{ab}

Diagram~1 in the Introduction summarizes the ZFC results in Part I of
the paper. The independence results in Part II show that any two rows
in Diagram~1 can be separated.
\end{abstract}

\maketitle

\section {Introduction}

The infinite Ramsey theorem in its simplest form states that
whenever all unordered pairs from an infinite set $A$ are colored
by two colors, there exists an  \emph{infinite homogeneous} $B\su
A$: an infinite  subset $B\su A$ with all  unordered pairs from it
colored by the same color.   Sierpinski constructed
pair-colorings on $\R$ with respect to which every homogeneous
set is countable, thus showing that there is no better Ramsey
theorem on $\R$ than there is on $\mathbb N$.

It is not too hard to check that if one colors all pairs from the
continuum  by two colors \emph{continuously} with respect to some
 complete, separable metric topology, then there is always a
nonempty perfect, hence  of size continuum, homogeneous set, and
that, furthermore, the \emph{chromatic number} of the coloring is
either countable or $2^{\aleph_0}$.

This fact shows that a Ramsey theorem on the continuum holds  for
continuous colorings,
  but also implies that for such colorings the standard Ramsey
  invariants \emph{clique number} and \emph{chromatic number} are degenerate,
from a set-theoretic point of view, being either countable or
equal to the continuum. (This  holds  also for open colorings on
analytic sets \cite{feng}.)

Recently a third  Ramsey invariant of continuous
colorings appeared
 in the classification of convex covers of closed planar
sets. For some closed subsets of $\R^2$ the number of convex
subsets required to cover them is equal to the \emph{homogeneity
number} $\hm(c)$ of some continuous pair-coloring $c$ on the Baire
space \cite{GKKS}. The homogeneity number is the least number of
homogeneous sets (of both colors) required to cover the space.

Unlike the chromatic  and clique numbers, homogeneity numbers of
continuous pair-colorings on the continuum  are \emph{not}
set-theoretically degenerate. Their classification leads to an
interesting theory in ZFC and to  two new forcing notions.

The broader class of \emph{open colorings} has been a focus of
interest for set theorists for three decades now,  and motivated
several important developments in the technique of forcing
\cite{B1,B2,ARS}. \emph{Open coloring axioms}, which are
statements in the Ramsey theory of open colorings, are among the
more frequently used set-theoretic axioms in the theory of the
continuum (see \cite{boban,Todorcevic,FA,JTM} and the references therein).

The crucial inequality (Theorem \ref{geqd} below) which enables
the reduction of the classification of  general continuous
pair-colorings by reducing them to  compact ones  involves the
notion of \emph{covering a plane by functions}.  About half of
the paper is devoted to that subject. The connection between
continuous pair-colorings and covering a plane by functions works
in both ways:  after establishing the classifications of
homogeneity numbers we have at hand an \emph{optimal forcing} for
proving the consistency of
 ``more Lipschitz functions are required to cover $\R^2$ than
continuous ones''.

\subsection{The results}
Two simple pair-colorings $c_{\min}$ and $c_{\max}$  are defined
on the Cantor space, and are shown to satisfy for every Polish space~$X$
and every continuous $c:[X]^2\to 2$ with uncountable $\hm(c)$:

\begin{equation}\label{minmaxint}
\hm(c_{\min})\le \hm(c)\le \hm(c_{\max})
\end{equation}

To state the remaining results concisely, we briefly  introduce
some notation. A function $f:X\to X$ \emph{covers} a point
$(x,y)\in X^2$ if $f(x)=y$ or $f(y)=x$. For a metric space
$(X,\dist)$ let
$\cov(\Lip(X))$ denote the number of Lipschitz functions from $X$
to $X$ required to cover $X^2$  and $\cov(\Cont(X))$ denote the
analogous  numbers for continuous functions. The \emph{Baire
space} $\om^\om$ and the \emph{Cantor space} $2^\om$ are
considered with the standard metric $\dist(x,y)= \frac {1}
{2^{\Delta(x,y)}}$, where
$\Delta(x,y)=\min\{n:x(n)\not=y(n)\}$ for $x\not=y$.

\medskip

The remaining ZFC equalities and inequalities are summarized in
Diagram~1.

\begin{figure}
$$
\begin{array}{cccc}
 (6) & &&\hsqueeze{ \bigg( \cov(\Cont(2^\omega))\bigg)^ +}\\
     & &&   \LineNotArrow\\
 (5) & && 2^{\aleph_0}\\
     & &&   \LineNotArrow\\
 (4) & && \hm(c_{\max})\\
     & &&   \LineNotArrow\\
 (3) & \nop{\cov(\Lip({\mathbb R})) \ge \cov(\Lip(\om^\om)) =
\cov(\Lip(2^\omega))} & =
& \hm(c_{\min})\\
     &   \LineNotArrow \\
 (2)    &
\nop{\cov(\Cont(2^\omega))=\cov(\Cont(\omega^\omega))=\cov(\Cont(\mathbb R))}
\\
     &   \LineNotArrow\\
 (1) & \d\\
 %      & x\qquad \qquad \qquad \qquad \qquad \qquad y
     & \end{array}
$$
\vskip0.5cm
\begin{center} Diagram 1
\end{center}

\end{figure}
\bigskip Homogeneity numbers are on the right column and
covering-by-functions cardinals are on the middle column. We draw
attention to the fact that the rows $(2)$--$(6)$ have to share at
most two consecutive cardinals since $\cov(\Cont(2^\om))$ cannot be more
than one cardinal below $2^{\aleph_0}$; thus, four different models of
set theory are required to separate them from each other.

The independence results in  Part II of the paper show that for
each of the rows $(1)$--$(5)$  it is consistent that the value at
the row is $\aleph_1$ and at all rows above the value is $\aleph_2$. The
forcing for separating $(2)$ from $(3)$ is a new example of an
\emph{optimal forcing} in the sense of Zapletal \cite{zapletal}
for increasing a cardinal invariant while leaving small
everything that can be left small.

\medskip
The inequality $\cov(\Lip(2^\omega))\le\hm(c_{\min})$ and the
consistency of $\hm(c)<2^{\aleph_0}$ for every Polish space $X$ and
 continuous $c:[X]^2\to 2$ were proved in \cite{GKKS}.

\medskip
 The last inequality cannot hold for all
\emph{open} colorings.   In \cite{ARS} an
 example of an open pair-coloring on  the square of any
 uncountable Polish space $X$
 is given such that $X^2$ cannot be covered by fewer than
 $2^{\aleph_0}$ homogeneous sets.  Let us present a slightly simplified
 version of
 this coloring.

 An unordered pair $\{(x_0,y_0),(x_1,y_1)\}$ of elements of $X^2$
 is of color $0$ if it is a 1-1-function and  of color $1$ otherwise.
 The set of pairs of color $0$ is open.  If $H\subseteq X^2$ is
 homogeneous of color~$1$, then it is either (a part of) a row or (a part
 of) a column in the square.  The homogeneous sets of color $0$ are graphs
 of (partial) injective functions.  It is easily checked that
 $X^2$ cannot be covered by less than $2^{\aleph_0}$ homogeneous
 sets.

\medskip

\subsubsection{Structure of the paper}

The paper is divided to two parts. Absolute ZFC results are in
Part I and  independence results are in Part II. Notation,
preliminaries and background material are included at the
beginning of each section. The first part employs elementary
techniques and does not require any specialized knowledge.

Although we are supposed to assume that every reader will read
the whole paper, we suspect that those who will read the second
part are knowledgeable in forcing notation. For those readers who
read the first part and decide that they have to learn forcing so
that they can read the second part, we recommend the standard
\cite{Kunen, Bartoszynski} as sources for notation and
introduction to forcing.

We tried to keep
notation as standard as possible.

\newpage

\begin{center}

\Large{\sc Part I: Results in ZFC}

\end{center}
\section{The structure of Continuous pair-colorings on Polish spaces}

\subsection{Basic definitions and preliminary facts}

\subsubsection{Colorings, chromatic numbers, and homogeneity numbers}
The symbol $[A]^2$ denotes the set of all two-element subsets of
a set~$A$.  Ramsey's theorem states that if $A$ is  infinite, then
for every function $c:[A]^2\to 2:=\{0,1\}$ there is an infinite
set $B\su A$ so that $c$ is constant on $[B]^2$. A function
$c:[A]^2\to 2$ is called a \emph{pair-coloring}, and a set $B\su
A$ for which $c\restriction[B]^2$ is constant is called
\emph{$c$-homogeneous} or \emph{$c$-monochromatic}. In the future
we write just $c\restriction B$ instead of $c\restriction[B]^2$.
A set $H$ is \emph{$c$-homogeneous of color $i$} for $i\in 2$, if the
constant color on $H$ is~$i$.

A pair coloring $c$ on $A$ can be thought of as (the
characteristic function of) the edge relation of a graph
$G=(A,c)$. In this setting  Ramsey's theorem states that every
infinite graph  contains either an infinite \emph{clique}
--- a subgraph in which any pair of vertices forms an edge ---
or  an infinite \emph{independent set} --- a subset in which no
two vertices form an edge.

Recall that the \emph{chromatic number} of a graph is the least
number of independent sets required to cover the set of vertices.

\begin{defn} For a coloring $c:[A]^2\to 2$ the \emph{homogeneity
number of $c$}, denoted by $\hm(c)$, is the minimal number of
$c$-homogeneous subsets required to cover~$A$.
\end{defn}

The difference between chromatic  and homogeneity numbers is that
in the definition of the latter  covering is  by homogeneous sets
of both colors.

\subsubsection{Continuous colorings on Polish spaces} Let $X$ be a
topological space and
let $X^2:=X\times X$ with the product topology.
We identify $[X]^2$ with
the quotient space $(X^2\setminus \{(x,x):x\in X\})/\sim$,
where $(x,y)\sim (w,z)$ iff $(x,y)=(w,z)$ or $(x,y)=(z,w)$.

A coloring $c:[X]^2\to 2$ is continuous if
the
preimages of $0$ and of $1$ are open.  Equivalently, $c$ is
continuous if for all $\{x,y\}\in[X]^2$, there are disjoint open
neighborhoods $U$ and $V$ of $x$ and~$y$, respectively, such that
$c$ is constant on $U\times V$.  Here we identify $c$ with the
corresponding
symmetric function from $X^2\setminus\{(x,x):x\in X\}$ to~$2$.

A topological space $X$ is \emph{Polish} if it is homeomorphic to a
separable and complete metric space.  Every Polish space is a
disjoint union of a countable open \emph{scattered} subset with a
\emph{perfect} subset (where either of the two components may be empty).
 Since every nonempty perfect subset of a
Polish space has the cardinality of the continuum, every
uncountable Polish space is equinumerous with the continuum.

\begin{defn}
\begin{enumerate}
\item  A  pair-coloring $c$ on $X$ is {\em reduced} if
$c$ is continuous and no nonempty open subset of $X$ is
$c$-homogeneous.
\item A coloring $c:[X]^2\to 2$ is \emph{trivial} if $\hm(c)\le
\aleph_0$.
\end{enumerate}
\end{defn}

\begin{fact} \label{reducedfact} If $X$ is a Polish space and
$c:[X]^2\to 2$, then $X=X_0\cup X_1$ such that $X_0$ is open,
$X_1$ is perfect, $X_0\cap X_1=\emptyset$, $c\restriction X_0$ is
trivial and $c\restriction X_1$ is reduced.
\end{fact}

\begin{proof} Let
$X_0$ be the union of all open sets $U\su X$ for
which $c\rest U$ is trivial. $X_0$ is open and since $X$ has a
countable basis, $c$ is trivial on $X_0$. Let $X_1=X\sm X_0$.
\end{proof}

\begin{fact} A continuous pair-coloring on a Polish space $X$
satisfies $\hm(c)>\aleph_0$ if and only if there exists a nonempty perfect
$Y\su X$ so that  $\hm(c)=\hm(c\rest Y)$ and $c\rest Y$ is reduced.
\end{fact}

\begin{proof} Suppose $\hm(c)>\aleph_0$ and write $X=X_0\cup X_1$ as stated
in the previous Fact. So $c\rest X_1$ is reduced.  Since $\hm(c\rest
X_0)\le \aleph_0$ it follows that $\hm(c)=\hm(c\rest X_1)$
and clearly $X_1\not=\emptyset$. On the other hand, suppose $Y\su X$
is perfect and nonempty, that $\hm(c)=\hm(c\rest Y)$ and $c\rest
Y$ is reduced.  Continuity of $c$ gives that the closure of every
$c$-homogeneous set is again $c$-homogeneous; so if  $Y\su X$ is perfect
and
$c\rest Y$ is reduced, every $c$-homogeneous subset of $Y$ is nowhere
dense and by the Baire theorem $\hm(c)>\aleph_0$.
\end{proof}

\subsubsection{Notation} Let  $\omega^\omega$ denote the set of all
(infinite) sequences of natural numbers. Let $\omega^{<\om}$
denote the set of all finite sequences of natural numbers and let
$\omega^{\le \omega}=\om^{<\om}\cup \om^\om$. Similarly,
$2^\omega,2^{<\omega}, 2^{\le \omega}$ are the analogous sets for
sequences over $\{0,1\}$.

\begin{defn}For $x,y\in\omega^{\leq\omega}$ let
$\Delta(x,y)=\min\{n\in\omega:x(n)\not=y(n)\}$ if there is some
$n\in\omega$ such that $x(n)\not=y(n)$.  Otherwise $\Delta(x,y)$
is undefined.

 If $\Delta(x,y)$ is defined for $x,y\in \om^{\le \omega}$, put
 \[\dist(x,y):=\frac 1 {2^{\Delta(x,y)}}\]
If $\Delta(x,y)$ is not defined, put $\dist(x,y):=0$.
\end{defn}

The function $\dist$ satisfies the triangle inequality. In fact,
it satisfies a stronger inequality: $\dist(x,z) \le
\max\{\dist(x,y),\dist(y,z)\}$ for all $x,y,z$.  (This makes $\dist$ an
\emph{ultra-metric}.)

The following  Polish spaces play an important role in this
section: the \emph{Cantor space} $(2^\om,\dist)$ and the
\emph{Baire space} $(\om^\om,\dist)$. These spaces are indeed
complete, separable metric spaces. The Cantor space is
homeomorphic to the usual \emph{Cantor set} and the Baire space
is homeomorphic the the space of \emph{irrational numbers}.

\subsubsection{The minimal coloring $c_{\min}$}

\begin{defn} If $X$ and $Y$ are topological spaces and $c$ and $d$ are
continuous pair-colorings on $X$ and~$Y$, respectively, then we
write $c\leq d$ if there is a topological embedding $e:X\to Y$,
such that for all $\{x_0,x_1\}\in[X]^2$,
$c(x_0,x_1)=d(e(x_0),e(x_1))$.
\end{defn}

  Clearly, if
$c\leq d$    via an embedding $e:X\to Y$, then $e^{-1}[A]$ is
$c$-homogeneous for every $d$-homogeneous $A\su Y$. Hence, $c\le
d$ implies that $\hm(c)\leq\hm(d)$.

We introduce next a pair coloring $c_{\min}$ on the Cantor space
which satisfies  $c_{\min}\le c$ for all reduced~$c$.

\begin{defn}\label{defcmin}
\begin{enumerate}
\item Let\/ $\parity(x,y)$ denote the parity of\/ $\Delta(x,y)$ for
$x,y\in \omega^{\le \omega}$ such that  $\Delta(x,y)$ is defined.

\item Let $c_{\parity}:=\parity\restriction\omega^\omega$.
\item Let
$c_{\min}:=\parity\restriction 2^\omega$.
\end{enumerate}
\end{defn}

Clearly, $c_{\parity}$ is a reduced pair-coloring on
$\omega^\omega$ and $c_{\min}$ is a reduced pair-coloring  on
$2^\omega$.

If $H\su 2^\omega$ is $c_{\min}$-homogeneous of color 0, then all
splittings in~$T(H)$, the tree of all finite initial segments of
members of~$H$, occur on even levels.  If $T$ is a subtree of
$\omega^{<\omega}$, we identify every infinite branch of $T$ with its
union, a point in $\omega^\omega$.   A set $H\su 2^\om$ is,
then, maximal $c_{\min}$-homogeneous  of color 0 is if and only
if  $H$ is the set of all infinite branches of a tree $T$ in
which $t\in T$ has two immediate successors if $|t|$ is even and
one immediate successor if $|t|$ is odd. Similarly, $H$ is
maximal $c_{\min}$-homogeneous of color 1 if and only if it is
the set of all infinite branches of a tree $T$ such that $t\in T$ has two
immediate successors in $T$ if $|t|$ is odd and one immediate
successor in $T$ if $|t|$ is even.

\begin{lemma}\label{canmin}  For every reduced pair-coloring $c$ on a
Polish space we have: $$c_{\min}\leq c.$$
 Consequently,
$\hm(c_{\min})\le \hm(c)$ for every reduced~$c$.
\end{lemma}

\begin{proof}
Suppose $c:[X]^2\to 2$ is reduced and $X$ is Polish. Since no
nonempty open set is $c$-homogeneous in~$X$, $X$ has no isolated
points.

By induction on $n$ choose, for every $t\in 2^n$, an open set
$U_t\not=\emptyset$ of diameter $<1/n$ such that
\begin{enumerate}
\item[--]  $t\su s\Rightarrow \cl(U_s)\su U_t$,
\item[--] $\Delta(t_1,t_2)$ defined implies
that $\cl(U_{t_1})\cap \cl(U_{t_2})=\emptyset$, and
\item[--]  for
every $t_1,t_2$ and $x_1\in \cl(U_{t_1})$, $x_2\in \cl(U_{t_2})$:
$c(x_1,x_2)\equiv n \mod 2$.
\end{enumerate}
At the induction step,
for a given $t\in 2^n$ find $x_1,x_2\in U_t$ which satisfy
$c(x_1,x_2)\equiv n \mod 2$ (possible since $U_t$ is not
$c$-homogeneous) and inflate $x_1,x_2$ to a sufficiently small
open balls $U_{t\concatenate 0},U_{t\concatenate 1}$.

The map $e$ mapping each $x\in 2^\om$ to the unique element of
$\bigcap_n U_{x\rest n}$ is an embedding of $2^\om$ into $X$
which preserves $c_{\min}$.
\end{proof}

In \cite{GKKS} $\hm(c_{\min})$ was denoted simply by~$\hm$.
We will also sometimes write $\hm$ for $\hm(c_{\min})$.

Before we proceed, let us remark   that $c_{\parity}$ is not more
complicated than $c_{\min}$:

\begin{lemma}  $c_{\parity}\leq c_{\min}$
\end{lemma}

\begin{proof} We have to define an embedding $e:\omega^\omega\to 2^\omega$
witnessing $c_{\parity}\leq c_{\min}$.

For $x\in\omega^\omega$, let $e(x)$ be the concatenation of the
sequences~$b_n$, $n\in\omega$, which are defined as follows.

If $n$ is even, then let $b_n$ be the sequence of length $2\cdot
x(n)+2$ which starts with $2\cdot x(n)$ zeros and then ends with
two ones.  If $n$ is odd, let $b_n$ be the sequence of length
$2\cdot x(n)+2$ starting with $2\cdot x(n)+1$ zeros and ending
with a single one.

It is clear that $e$ is continuous and it is easy to check that
$e$ is an embedding witnessing $c_{\parity}\leq c_{\min}$.
\end{proof}

\subsection{Classification of homogeneity numbers} We begin now
the classification of
homogeneity numbers of continuous pair colorings on
Polish spaces. The following
sequence of reductions will be performed: From general Polish
spaces to compact metric spaces; from compact metric spaces to
the Cantor space; and from the class of all continuous pair
colorings on the Cantor space to a subclass of particularly
simple colorings.

\subsubsection{Reduction to compact spaces}
The following two fundamental inequalities hold for $c_{\min}$:

\begin{gather}
(\hm(c_{\min}))^+\ge 2^{\aleph_0}\label{one}\\
\hm(c_{\min})\ge \mathfrak d \label{two}
\end{gather}

The first inequality was proved in \cite{GKKS} and the second one
which, really, is the starting point of the present paper, will
be proved in  Section \ref{covsection}. Although these
inequalities are central for this Section, their
proofs belong to the setting of covering a square by functions.

 From the first inequality it
follows that there is room for at most one more homogeneity
number above $\hm(c_{\min})$ --- since either $\hm(c_{\min})$ or
its immediate successor cardinal is the continuum. In \cite{GKKS}
it was proved consistent that for all reduced pair-colorings~$c$,
\begin{equation}
\hm(c)=\aleph_1<2^{\aleph_0}=\aleph_2.
\end{equation}

The second inequality relates $\hm(c_{\min})$ to the
\emph{domination number} $\mathfrak d$. This number is the least
number of functions from $\omega$ to $\omega$ needed to
eventually dominate every such function. Another important
feature of $\mathfrak d$ is that $\omega^\omega$ can be covered
by $\mathfrak d$ compact sets. It is well-known that every Polish
space is a continuous image of $\omega^\omega$.    Therefore
every Polish space can be covered by $\mathfrak d$ compact sets.

\begin{lemma}\label{redcomp}
For every Polish space $X$ and a continuous pair-coloring $c:[X]^2\to
2$ with uncountable $\hm(c)$ there is a compact subspace $Y\su X$ so
that $\hm(c)=\hm(c\rest Y)$.
\end{lemma}

\begin{proof}
 Suppose, without loss of generality, that $c$ is reduced on
$X$. Cover $X$ by compact subspaces $Y_\alpha$, $\alpha\le \mathfrak
d$, and denote $c_\alpha:=c\rest Y_\alpha$. For each $\alpha<\mathfrak
d$ fix a collection $\mathcal U_\alpha$ of $c_\alpha$-homogeneous
subsets of $Y_\alpha$ which covers $Y_\alpha$ and such that $|\mathcal
U_\alpha|=\hm(c_\alpha)$.  Thus $\mathcal U =\bigcup_{\alpha<\mathfrak
d}\mathcal U_\alpha$ is a collection of $c$-homogeneous sets which
covers~$X$, so $ \hm(c)\le |\mathcal U|$.

In the case that for all $\alpha<\mathfrak d$ it holds that
$\hm(c_\alpha)\le\hm(c_{\min})$ we have  that  $\hm(c)
\le|\mathcal U|\le \mathfrak d\cdot\hm(c_{\min})$, so by
(\ref{two}), $\hm(c)\le \hm(c_{\min})$. Since $c$ is reduced,
$\hm(c)=\hm(c_{\min})$ and $Y\su X$ can be chosen as a copy of
the Cantor space by Lemma \ref{canmin}

In the remaining case $\hm(c)>\hm(c_{\min})$,  therefore there
necessarily exists $\alpha<\mathfrak d$ for which
$\hm(c_\alpha)>\hm(c_{\min})$, and consequently, by (\ref{one}),
$\hm(c_\alpha)=\hm(c)$.
\end{proof}

 Now it is clear, subject to  the inequalities above,
that all homogeneity numbers of continuous pair-colorings on
arbitrary Polish spaces appear on compact Polish spaces, i.e., on
compact metric spaces.

\subsection{Pair-colorings on compact metric spaces}

In this Section we reduce the study of  continuous a
pair-colorings on compact metric spaces to continuous
pair-colorings on $2^\omega$, and then reduce it further to a
class of particularly simple colorings on $2^{\omega}$. At the
end of the section, we shall be able to isolate a pair-coloring
$c_{\max}$ with a maximal homogeneity number in the class of
continuous pair-colorings on Polish spaces.

\subsubsection{Getting rid of topological connectedness}
For a compact space let $\Comp(X)$ be the set of connected
components of~$X$.  For $x\in X$ let $\comp(x,X)$ denote the
component of $x$ in $X$ and $\comp(x)=\comp(x,X)$ when $X$ is
clear from the context. $\Comp(X)$ becomes a compact space when
equipped with the quotient topology.

The components of $\Comp(X)$  are singletons.  Since $\Comp(X)$ is
compact, it is zero-dimensional.  (See
\cite{engelking} for this.)

\begin{lemma}\label{zerodim} Let $X$ be compact and $c:[X]^2\to 2$
continuous. Define a coloring $\overline c:[\Comp(X)]^2\to 2$ by
\[ \overline c\bigl(\comp(x),\comp(y)\bigr)=c(x,y)\]
for all $x,y\in X$ with $\comp(x)\not=\comp(y)$. Then $\overline
c$ is a well-defined  continuous pair-coloring on
$\Comp(X)$.
\end{lemma}

\begin{proof}
Suppose $x_0,x_1,y_0,y_1\in X$ are such that
$x_1\in\comp(x_0)$, $y_1\in\comp(y_0)$, and
$x_0$ and $y_0$ are in different components.  Then
$c(x_0,y_0)=c(x_1,y_0)$ since $x_0$ and $x_1$ are in the same
component of $X\setminus\{y_0\}$ and
$c(\cdot,y_0):X\setminus\{y_0\}\toÝ2$ is continuous.  By the same
argument, $c(x_1,y_0)=c(x_1,y_1)$. Thus $c(x_0,y_0)=c(x_1,y_1)$,
showing that $\overline c$ is well-defined.

For every $x\in \comp(x_0), y\in\comp(y_0)$ fix, by continuity of
$c$, disjoint open $U_{x,y}\ni x$, $V_{x,y}\ni y$ so that $c$ is
constant on $U_{x,y}\times V_{x,y}$. Now $\{U_{x,y}\times
V_{x,y}:x\in \comp(x_0), y\in\comp(y_0)\}$  is an open cover of
$\comp(x_0)\times \comp (y_0)$. Since the latter is compact,
there is a finite subcover $\{U_{x_i,y_i}\times
V_{x_i,y_i}:i<n\}$ of this cover, which can be shrunk so that
$\bigcup_{i<n} U_{x_i,y_i}\cap \bigcup _{i<n}
V_{x_i,y_i}=\emptyset$. Thus we found two disjoint open neighborhoods of
$\comp(x)$, $\comp(y)$ respectively so that $c$ is constant on
their product. This proves the continuity of $\overline c$.
\end{proof}

Recall that in a compact space the connected component of a point
is equal to the intersection of all clopen sets that contain the
point (see \cite{engelking}).

\begin{lemma}\label{connectedsets} Let $X$ be compact and connected.  Then
every continuous $c:[X]^2\to 2$ is constant.  In other words,
$[X]^2$ is connected.
\end{lemma}

\begin{proof} We need:

\begin{claim} Suppose $X$ is compact and connected and let $x\in
X$. Then for every $y\in X\sm \{x\}$, the point  $x$ is in the
closure (in $X$) of  $\comp(y,X\sm\{x\})$.
\end{claim}

\begin{proof}
Let $y\in X\sm \{x\}$ be arbitrary and let $Y:=\comp(y,X\sm
\{x\})$. If $x\notin\cl_X(Y)$, then $Y$ is closed in~$X$. By
normality of~$X$, there is an open $U\ni x$ (in $X$) so that $\cl_X
(U)\cap Y=\emptyset$. Replacing $U$ by $\int\cl_X(U)$ we may assume
that $U$ is regular open, therefore $\bd_X(U)=\cl_X(U)\sm U$.

The space $X\sm U$ is compact, and $Y$ is the
component of $y$ also in $X\sm U$.  The sets $Y$ and $\bd_X(U)$
are closed and disjoint subsets of $X\sm U$, so since $Y$ is an
intersection of clopen sets, there is, by compactness of $X\sm U$, a
\emph{finite} intersection $V$ of clopen sets, thus itself clopen,
which contains $Y$ and is disjoint from $\cl_X(U)$. Thus $V$ is clopen
in $X$ and $X$ is not connected.
\end{proof}

 Suppose now that $c:[X]^2\to 2$ is not constant.
If $c$ is constant on the pairs from every 3-element subset of
$X$, it is constant; thus  there are distinct  $x,y,z\in X$ such
that $c(x,y)=0$ and $c(x,z)=1$.

Let $Y:=\comp(y,X\setminus\{x\})$ and $Z:=\comp(z,
X\setminus\{x\})$. By the previous claim, $x\in \cl_X(Y)\cap
\cl_X(Z)$.

 Since $Y$ is connected in $X\setminus\{x\}$ and does not
contain~$z$, $c(z,y')=c(z,y)$  for every $y'\in Y$. Since $x\in
\cl_X(Y)$, continuity of $c$ implies that $c(z,y)=c(z,x)=1$.
Symmetrically, $c(y,z)=c(y,x)=0$. Hence $0=c(y,z)=c(z,y)=1$
--- a contradiction.
\end{proof}

\begin{problem} Is it true for an arbitrary connected Hausdorff space
$X$ that
$[X]^2$ is connected?
\end{problem}

\subsubsection{Reduction to colorings on $2^{\omega}$}

\begin{lemma}\label{redtocantor} Let $X$ be a compact metric space and
suppose $c:[X]^2\to 2$ is continuous. Then there exists a
continuous $\overline c:[2^\omega]^2\to 2$ such that $\hm(c)\le
\hm(\overline c)$.
\end{lemma}

\begin{proof}
Let $Y:=\Comp(X)$ and let $f:X\to Y$ be the mapping that maps
every $x\in X$ to $\comp(x,X)$.  Let $\overline c$ be as in Lemma
\ref{zerodim}. Observe  that $Y$ is of countable weight.

Assume that $Y$ is uncountable. Cantor-Bendixson analysis of $Y$
gives us a decomposition of $Y$ into countably many points and a
perfect set.   Since for every isolated point $y\in Y$ the set
$f^{-1}(y)$ is $c$-homogeneous in $X$ by Lemma
\ref{connectedsets}, we may replace $Y$ by a perfect subset of $Y$
at the cost of removing countably many $c$-homogeneous subset of
$X$.

$Y$ is now zero-dimensional, compact, without isolated points and
of countable weight.  Therefore $Y$  is the Cantor space.

 \begin{claim} $\hm(c)\le \hm(\overline c)$
\end{claim}

By the continuity of $\overline c$, every maximal $\overline
c$-homogeneous set in $Y$ is closed. Now using Cantor-Bendixson
analysis again, every uncountable maximal $\overline
c$-homogeneous set can be decomposed into countably many
singletons and a perfect set.

The preimages under $f$ of singletons are $c$-homogeneous by Lemma
 \ref{connectedsets}. Also,

\begin{claim} For any perfect $\overline c$-homogeneous set $H\subseteq
Y$, $f^{-1}[H]$ is $c$-ho\-mo\-ge\-neous.
\end{claim}
\begin{proof}
For the claim let $H\subseteq Y$ be perfect and $\overline
c$-homogeneous of color $i\in 2$.  If $x,y\in f^{-1}[H]$ are in
different components of~$X$, then clearly $c(x,y)=i$.  Now let
$z$ be one of the components of~$X$. Assume $\card{z}>1$. By
Lemma \ref{connectedsets}, $c$ is constant on~$z$.  Let $j\in 2$
be the constant value of $c$ on~$z$. We have to show $i=j$.

Let $(z_n)_{n\in\omega}$ be a sequence in $H\setminus\{z\}$ that
converges to~$z$.  Pick $(x_n)_{n\in\omega}$ in $X$ such that for
all $n\in\omega$, $f(x_n)=z_n$. By compactness,
$(x_n)_{n\in\omega}$ has a convergent subsequence.  We may assume
that $(x_n)_{n\in\omega}$ itself converges.

Let $x$ be the limit of $(x_n)_{n\in\omega}$.  Clearly, $x\in
z$.  Let $y\in z$ be different from~$x$. Then $c(x,y)=j$.  By
continuity, $c(x,y)=\lim_{n\to\infty}c(x_n,y)=i$. Thus $i=j$,
which finishes the proof of the claim.
\end{proof}

Thus, the preimage under $f$ of every $\overline c$-homogeneous
subset of $Y$ is a countable union of $c$-homogeneous subsets of~$X$.
This establishes $\hm(c)\le\hm(\overline c)$ and proves the
theorem.
\end{proof}

\subsubsection{Reduction to simple colorings on $2^{\omega}$}
 We are now fishing in a much smaller tank: we can
consider only colorings on the Cantor space. The next reduction
will show that we can consider only ``coarse" pair-colorings on
the Cantor space.

\begin{Notation} \label{notation}  %% mg May 11
  For a tree $T$ and $t\in T$ let
$\succ_T(t)$ be the set of immediate successors of $t$ in~$T$.
Recall that if $A$ is a subset of $\om^\om$, then $T(A)$ denotes the set
of finite initial segments of the element of~$A$, a subtree of
$\om^{<\om}$.  If $T$ is a subtree of $\om^{<\om}$, then $[T]$
denotes the set of all elements of $\om^\om$ which have all their
finite initial segments in~$T$. $[T]$ is a closed subset of $\om^\om$.
In this way closed subsets of $\om^\om$ correspond to subtrees of
$\om^{<\om}$ without finite maximal branches.
\end{Notation}

A natural way to construct continuous pair-colorings on a subset
$A$ of $\omega^\omega$ is the following: To each $t\in T(A)$
assign a coloring $c_t:[\succ_{T(A)}(t)]^2\to 2$. Now for all
$\{x,y\}\in [A]^2$ let $t$ be the longest common initial segment
of $x$ and $y$ and put $c(x,y):=c_t(x\restriction
n+1,y\restriction n+1)$ where $n=\dom(t)$. Clearly, $c$ is
continuous.  We call a coloring which is defined in this way an
{\em almost node-coloring}.

 A  {\em node-coloring} on $A$  is obtained
by assigning a color to every node $t\in T(A)$ and then defining
the color of $\{x,y\}\in [A]^2$ to be the color of the longest
common initial segment of $x$ and~$y$. Equivalently, a
node-coloring is an almost node-coloring in which
$c_t:[\succ_{T(A)}(t)]^2\to 2$ is constant for all $t\in
T$.

Both
$c_{\min}$ and
$c_{\parity}$ are node-colorings.

Not every continuous pair-coloring on $\omega^\omega$ is an almost
node-coloring.  However, the following holds:

\begin{lemma} Let $c:[2^\omega]^2\to 2$ be
continuous. Then there is a topological embedding
$e:2^\omega\to\omega^\omega$ such that for every
$c_{\parity}$-homogeneous set $H\subseteq e[2^\omega]$, the
coloring $c^e\rest H$ which is induced on $H$ by $c$ via $e$ is
an almost node-coloring.
\end{lemma}

\begin{proof}
Let $n\in\omega$ and let $s,t\in 2^{n+1}$ be such that
$\Delta(s,t)=n$. Let $O_s$ and $O_t$ denote the basic open
subsets of $2^\omega$ determined by $s$ and~$t$, respectively.

Since $O_s\times O_t$ is compact and $c$ is continuous, there is
$m>n$ such that for all $(x,y)\in O_s\times O_t$, $c(x,y)$ only
depends on $x\restriction m$ and $y\restriction m$.

It follows that there is a function $f:\omega\to\omega$ such that
for all $\{x,y\}\in[2^\omega]^2$, $c(x,y)$ only depends on
$x\restriction f(\Delta(x,y))$ and $y\restriction
f(\Delta(x,y))$.  We can choose $f$ strictly increasing and such
that $f(0)\geq 1$.  For $n\in\omega$ let $g(n):=f^n(0)$.

Identifying $2^{<\omega}$ and $\omega$, we define the required
embedding $e:2^\omega\to\omega^\omega$ by letting
$e(x):=(x\restriction g(0), x\restriction g(1),\dots)$. Let
$E:=e[2^\omega]$. $c$ induces a continuous pair-coloring $c^e$ on
$E$ via~$e$. By the choice of~$f$, for $\{u,v\}\in[E]^2$,
$c^e(u,v)$ only depends on $u\restriction(\Delta(u,v)+2)$ and
$v\restriction(\Delta(u,v)+2)$. This is because if
$n=\Delta(u,v)$ and $x,y\in 2^\omega$ are such that $e(x)=u$ and
$e(y)=v$, then $\Delta(x,y)<g(n)$ and thus $c(x,y)$ only depends
on $x\restriction f(\Delta(x,y))$ and $y\restriction
f(\Delta(x,y))$.  But since $f$ is strictly increasing,
$f(\Delta(x,y))<f(g(n))=g(n+1)$.

Now let $H$ be a $c_{\parity}$-homogeneous subset of~$E$. The
$c_{\parity}$-homogeneity of $H$ implies that for all
$\{u,v\}\in[H]^2$, the restrictions of $u$ and $v$ to
$\Delta(u,v)+1$ uniquely determine the restrictions to
$\Delta(u,v)+2$.  Therefore, for all $\{u,v\}\in[H]^2$,
$c^e(u,v)$ only depends on $u\restriction(\Delta(u,v)+1)$ and
$v\restriction(\Delta(u,v)+1)$.

It follows that $c^e\restriction H$ is an almost node-coloring.
\end{proof}

\begin{corollary}\label{redalmostnode}
For every  continuous pair-coloring $c:[2^\omega]^2\to 2$, there
is an almost node-coloring $d$ on some compact subset of
$\omega^\omega$ such that $\hm(c)\leq\hm(d)$.
\end{corollary}

\begin{proof} By the previous Lemma, $2^\om$ can be presented as a
union of $\le \hm(c_{\min})$ sets on each of which $c$ is
reducible to an almost node-coloring. The rest of the proof is as
in the proof of Lemma \ref{redcomp}.
\end{proof}

\subsubsection{The coloring $c_{\max}$}
We shall now  define a maximal almost node-coloring.

Recall that the \emph{random} graph on $\om$ is, up to
isomorphism, the only homogeneous and universal graph in the
class of all graphs on~$\om$. (See \cite{ES} for some information
on the random graph.) Universality means: every graph $(\om,E)$ is
embeddable as an induced subgraph into the random
graph (in particular, every finite graph is embeddable as an
induced subgraph into a finite initial segment of the random
graph).

\begin{defn}
Let $\chi_{\random}:[\omega]^2\to 2$ be the (characteristic
function of the) edge relation of the random graph.
For $s,t\in \om^{\le \omega}$ write $\random (s,t)=i$ iff
$n:=\Delta(s,t)$ exists and $i=\chi_{\random}(s(n+1), t(n+1))$.
Let $c_{\random}:[\omega^\omega]^2\to 2$ be defined by
$c_{\random}(x,y):=\random(x,y)$. Finally, let
\begin{equation}\label{cmax}
c_{\max}:=c_{\random}\restriction\prod_{n\in\omega}(n+1)
\end{equation}
\end{defn}

Clearly, $c_{\random}$ and $c_{\max}$ are almost node-colorings.
Since $\prod_{n\in\omega}(n+1)$ is homeomorphic to $2^\omega$, we
regard $c_{\max}$ as a coloring on $2^\omega$.

It is interesting to point out:

\begin{fact}
Whenever $c$ is an almost node-coloring
on a compact subspace of $\om^\om$, then:  $c_{\random}\not \le c$.
\end{fact}

\begin{proof} Let $(x_n)_{n\in \om}$ be an infinite path in
$c_{\random}$, i.e.,
$$ \forall n<m:\qquad   c_{\random}(x_n,x_{m})=1
    \  \Leftrightarrow\ m=n+1.$$
 Since every countable graph embeds into
$(\om^\om,c_{\random})$, such a sequence can be easily found.

On the other hand, if $Y\su \om^\om$ is compact and $c:[Y]^2\to 2$
is an almost node-coloring, there is no infinite path in $(Y,c)$. Suppose
to the contrary that $(y_n)_{n<\om}$ is a
path in $(Y,c)$. Observe that
$\Delta(y_{n+1},y_{n+2})>\Delta(y_{n},y_{n+1})$ implies that
$c(y_{n},y_{n+2})=1$; and that $\Delta(y_{n+1},y_{n+2})<
\Delta(y_{n},y_{n+1})$ implies $c(y_{n+1},y_{n+2})=0$. Thus,
$\Delta(y_n,y_{n+1})$ is constant for all $n$ --- contrary to the
compactness of~$Y$.

The fact now follows.
\end{proof}

\begin{lemma}\label{maxmax} a) If $c$ is an almost node-coloring on a subset of
$\omega^\omega$, then $c\leq c_{\random}$ via a level preserving
embedding (isometry) of $\om^\om$ into $\om^\om$.

b) If $c$ is an almost node-coloring on a compact subset of
$\omega^\omega$, then $c\leq c_{\max}$.
\end{lemma}

\begin{proof}
Let us prove b) first. Suppose $c$ is an almost node-coloring on
a compact subset $A$ of $\omega^\omega$.  Then $T(A)$ is a
finitely branching subtree of $\omega^{<\omega}$.  For each $t\in
T(A)$ fix a coloring $c_t:[\succ_{T(A)}(t)]^2\to 2$ such that the
$c_t$ witnesses the fact that $c$ is an almost node-coloring.
For $s,t\in T$ let $\overline{c}(s,t):=c(x,y)$ if $s$ and $t$ are
incomparable and $x,y\in[T]$ are such that $s\su x$ and $t\su
y$.  If $s$ and $t$ are comparable, then $\overline{c}(s,t)$ is
undefined.

Let  $T_k=\{t\in T(A): |t|=k\}$. We construct a monotone (i.e.,
$\su$-preserving) map $e:\bigcup_{k\in\om}T_k\to
T(\prod_{n\in\om}(n+1))$ which induces the  required embedding of
$A$ into $\prod_{n\in\om}(n+1)$.

Argue by induction on~$k$. Suppose that $e(s)\in\prod_{n\leq
n(k)}(n+1)$ is defined for all $s\in T_k$, and for all $s,t\in
T_k$ we already have  $\random(e(s),e(t))=\overline{c}(s,t)$. \\
Find
$n(k+1)>n(k)$ such that for all $s\in T_k$ there is
$t\in\prod_{n<n(k+1)}(n+1)$ with  $e(s)\su t$ and
$c_s\leq\random\rest\succ_{T(\prod_{n\in\omega}(n+1))}(t)$. Now it
is obvious how to define $e$ on $T_{k+1}$ with images in
$\prod_{n\leq n(k+1)}(n+1)$.

a) is proved similarly, using the fact that every countable graph
occurs as an induced subgraph of $(\Succ_{\om^{<\om}}(s),
\random)$ for every $s\in \om^{<\om}$.
\end{proof}

\begin{corollary}
For every Polish $X$ and every continuous $c:[X]^2\to 2$:
$$\hm(c)\le \hm(c_{\max}).$$
\end{corollary}

\begin{proof} Let $c$ be an arbitrary reduced continuous
pair-coloring on a Polish~$X$. By Lemma \ref{redcomp} there exists
a compact $Y\su X$ so that $\hm(c)=\hm(c\rest Y)$. By Lemma
\ref{redtocantor} there is a coloring $\overline c$ on $2^\om$ so
that $\hm(c)\le \hm(\overline c)$  and by Corollary
\ref{redalmostnode} there is an almost node-coloring $d$ on
$2^\om$ so that $\hm(\overline c)\le \hm(d)$. Finally, $d\le
c_{\max}$ by  Lemma \ref{maxmax} above.
\end{proof}

Finally,

\begin{thm}\label{minormax}
For every reduced continuous pair-coloring~$c$:
\[ \hm(c)=\hm(c_{\min}) \text{  or  } \hm(c)=\hm(c_{\max}) \]
\end{thm}

\begin{proof}
By now we have that $\hm(c_{\min})\le \hm(c)\le\hm(_{\max})$ for
all reduced~$c$. But $\hm(c_{\max})\le (\hm(c_{\min}))^+$ by
(\ref{one}); so $\hm(c)>\hm(c_{\min})$ implies
$\hm(c)=\hm(c_{\max})$.
\end{proof}

We remark that in Theorem \ref{minormax} above, $c_{\min}$ can be
replaced by $c_{\parity}$ and $c_{\max}$ can be replaced by
$c_{\random}$, since

\[c_{\parity}\le c_{\min} \le  c_{\max} \le c_{\random}. \]

\subsubsection{Why  $c_{\max}$ is more complicated than
$c_{\min}$: Random versus perfect graphs.}
 In the second part of the paper we
shall prove the consistency of $\hm(c_{\min})<\hm(c_{\max})$. The
consistency proof relies on the different   finite patterns that
appear in each of those two colorings.

 Clearly, every finite graph occurs
as an induced subgraph of $(2^\omega,c_{\max})$.

A finite graph is called \emph{perfect} if  in each of its induced
subgraphs the chromatic number is equal to the clique number. A
perfect graph with $n$ vertices contains either a clique or an
independent set of size $\lfloor \sqrt n\rfloor$. This stands in
strong contrast to a randomly chosen graph: in a random graph on
$n$ vertices there is almost certainly no clique and no
independent set of size $2\log n$ (see \cite{AS}).

\begin{fact}[N. Alon]
Every finite (induced) subgraph $H$ of $(\om^\om,c_{\parity})$
satisfies that the chromatic number of
$H$ is equal to the maximal size of a clique in~$H$.
\end{fact}

\begin{proof} Two proofs of this fact are in \cite{Alon}.
The proof we include here was suggested to us by Stevo Todor\v{c}evi\'c.
Define a partial order on $\om^\om$ by $\eta_1\le \eta_2$ iff
$\eta_1=\eta_2$ or $\Delta(\eta_1,\eta_2)$ is odd and $\eta_1$
precedes $\eta_2$ in the lexicographic ordering on $\om^\om$. A finite
induced subgraph of $\om^\om$ is a clique iff its elements form a
chain in the poset just defined and is an independent set iff its
elements form an anti-chain in the same poset.  Now recall that a
finite partially ordered set with no chain of length $k+1$ is a union
of $k$ antichains.
\end{proof}

Thus  only perfect graphs occur as
finite induced subgraphs of $c_{\min}$.

In particular:
\begin{equation}
c_{\max}\not\leq c_{\min}.
\end{equation}

\section{Covering a square by functions}\label{covsection}

The problem of covering a Euclidean space by smaller geometric
objects is well investigated. Klee \cite{klee} proved that no
separable Banach space can be covered by fewer than $2^{\aleph_0}$
hyperplanes.  Stepr\=ans \cite{steprans} proved the consistency of
covering $\R^{n+1}$ by fewer than continuum \emph{smooth
manifolds of dimension $n$}.

We recall that a point $(x,y)\in X^2$ is \emph{covered} by a
function $f:X\to X$ if $f(x)=y$ or $f(y)=x$. By $f^{-1}$ we mean
the set $\{(y,x):f(x)=y\}$.  Thus $(x,y)$ is covered by $f$ iff
$(x,y)\in f\cup f^{-1}$. For a metric space $X$ denote by
$\cov(\Cont(X))$ the minimal number of continuous functions from
$X$ to $X$ needed to cover $X^2$ and by $\cov(\Lip(X))$ denote
the analogous number for Lipschitz functions.

Hart and van der Steeg showed the consistency of covering
$(2^\omega)^2$ by fewer than continuum continuous functions
\cite{hartsteeg}, a result that actually follows from Stepr\=ans'
result mentioned above using some easy arguments from the present
article. Ciesielski and Pawlikowski proved that $\mathbb R^2$ is
consistently  covered by fewer than continuum continuously
differentiable {\em partial} functions with perfect domains
\cite{CPA}.

In \cite{GKKS} it was shown that $(2^\omega)^2$ can consistently
be covered by fewer than continuum Lipschitz functions. Hart
asked whether $\cov(\Lip(2^\omega))$ can be different from
$\cov(\Cont(2^\omega))$.    Recently, Abraham and Geschke
\cite{turing} proved that it is consistent to cover $\R^{n+1}$ by
$\k$ $n$-ary continuous functions with $2^{\aleph_0}=\k^{+n}$.

Let us state the following folklore result that was brought to the
authors' attention by Ireneusz Rec\l aw (and which should be
well-known):

\begin{thm}\label{coveringbyfunctions} Let $\kappa$ be an infinite
cardinal.  Then the least number of functions from
$\kappa^+$ to $\kappa^+$ needed to cover $\kappa^+\times\kappa^+$ is
$\kappa$.
\end{thm}

\begin{proof}
For every $\alpha<\kappa^+$ fix a surjection $f_\alpha:\kappa\to
(\alpha+1)$. Now define, for $\beta<\k$,
$g_\beta(\alpha)=f_\alpha(\beta)$. The functions
$\{g_\beta:\beta<\kappa\}$ cover $\kappa^+\times \kappa^+$.

To show that $\kappa^+\times\kappa^+$ is not covered by less than
$\kappa$ functions, let $X$ be any infinite set and let $\mathcal F$ be a
family of functions on $X$ which covers~$X^2$.  Assume that
$\id_X\in\mathcal F$ and $\mathcal F$ is closed under composition of
functions.  For $x,y\in X$ let $x\leq_{\mathcal F}y$ iff there is
$f\in\mathcal F$ such that $f(y)=x$.

It is easily checked that $\leq_{\mathcal F}$ is a linear quasi-ordering.
For every $x\in X$ the set $\{y\in X:y\leq_{\mathcal F}X\}$ has size at
most $\card{\mathcal F}$.  It follows that $\card{X}$ is not greater than
$\card{\mathcal F}^+$.
\end{proof}

In  \cite{turing} a generalization of this to higher dimension is
proved.

This theorem implies that if the continuum is a successor
cardinal, then fewer than continuum functions suffice to cover
the square of the continuum.

\medskip
In the rest of this  section the connection between
$c_{\min}$-homogeneous sets and covering $(2^\om)^2$ by Lipschitz
functions will be explored, and used to prove the inequalities
(\ref{one}) and (\ref{two}) which were used in the previous
Section. Inequality (\ref{one}) was already proved in
\cite{GKKS}. Inequality (\ref{two}) follows from Theorem
\ref{geqd} below.

After proving the crucial  Theorem \ref{geqd} we investigate
covering by \emph{continuous} functions.

\subsection{$\hm$ and covering a square by Lipschitz functions}
\begin{defn}For $a,b>0$ let $\Lip_{a,b}$ denote the $\sigma$-ideal on
$2^\omega$ generated by the (graphs of) Lipschitz functions of
constant $a$ and the reflections on the diagonal (of graphs) of
Lipschitz functions of constant $b$ (i.e., {\em inverses} of Lipschitz
functions of constant $b$). The \emph{covering number} of
this ideal, $\cov(\Lip_{a,b})$, is the least number of sets in
the ideal needed to cover $(2^{\omega})^2$.
\end{defn}

Clearly, as the graph of every continuous function is a
nowhere-dense subset of $(2^\omega)^2$,
$\cov(\Lip_{a,b})>\aleph_0$ for every choice of positive~$a,b$.
By Theorem \ref{coveringbyfunctions} we know  that
$(\cov(\Lip_{a,b}))^+\ge 2^{\aleph_0}$.

\begin{lemma}\label{limitedlipeqhm} $\hm=\cov(\Lip_{1,\frac 1 2})$
\end{lemma}

\begin{proof}
 For $x,y\in\omega^\omega$ let $x\otimes
y:=(x(0),y(0),x(1),y(1),\dots)$.  It is easily seen that
$\otimes:(\om^\om)^2\to\om^\om$ and $\otimes:(2^{\om})^2 \to 2^\om$
are uniformly continuous homeomorphisms.

 Suppose $H_0\subseteq 2^\omega$ is a maximal $c_{\min}$-homogeneous
 of color
$0$. Then $T:=T(H_0)$ is a tree with the property that $t\in T$ has
two immediate successors in $T$ if and only if $|t|$ is even and
has one immediate successor in $T$ otherwise. Let $x\in 2^\omega$
and define $y(n)$ inductively as follows:\\
 Suppose $y(i)$ is
defined for all~$m<n$, and we have
$$t=
(x(0),y(0),x(1),y(1),\dots ,x(n-1),y(n-1)) \in T.$$
Let $y(n)\in
\{0,1\}$ be the unique
 such that $t \concatenate i\in T$. Let $f_{H_0}(x)$ denote~$y$,
which we have just defined from $x$ and~$H_0$. We then have
$(x\otimes f_{H_0}(x))\in H_0$.

 Since
the first $n$ digits of $y$ are determined by the first $n$ digits of
$x$,
$f_{H_0}:2^\om\to 2^\om$ is a Lipschitz function with constant 1
(with respect to $\dist$).

Similarly, if $H_1$ is maximal $c_{\min}$-homogeneous of color 1,
then for every $x\in 2^{\omega}$ there is a unique
$f_{H_1}(x)\in 2^{\omega}$ for which $f_{H_1}(x)\otimes x\in H_1$. This
time, the function $f_{H_1}$ is of Lipschitz of constant $\frac 1 2$.

Conversely, from every $1$-Lipschitz function $f:2^\omega\to
2^{\omega}$   a maximal $c_{\min}$-homogeneous set $H_f$ of color
$0$ is defined so that for all~$x$, $y=f(x)$ is the unique such
that $x\otimes y\in H_f$ and from every $1/2$-Lipschitz function
$f:2^\omega\to 2^{\omega}$ a maximal $c_{\min}$-homogeneous set
${}_f H$ of color $1$ is defined such that $y=f(x)$ is the unique
such that $y\otimes x\in {}_f H$.

Suppose $\mathcal H_0$ is a family of maximal
$c_{\min}$-homogeneous subsets of $2^\om$ of color $0$ and
$\mathcal H_1$ is a family of maximal $c_{\min}$-homogeneous
subsets of color 1.  For $(x,y)\in (2^\omega)^2$, if
$x\otimes y\in H$ for some $H\in \mathcal H_0$ then $y=f_{H}(x)$,
and if $x\otimes y\in H$ for $H\in \mathcal H_1$ then
$x=f_{H}(y)$. Thus  $\bigcup\mathcal H_0\cup \bigcup \mathcal
H_1=2^\omega$ implies that for all $(x,y)\in (2^\omega)^2$ there
is some $H\in \mathcal H$ for which $f_{H}(x)=y$ or $f_H(y)=x$.

Conversely, suppose that $\mathcal F_0$ is a family of 1-Lipschitz
functions from $2^\om$ to itself and that $\mathcal F_1$ is a
family of $\frac 1 2$-Lipschitz functions from $2^\om$ to itself. Let
$z\in 2^{\omega}$ and write
$z=x\otimes y$. If there is $f\in \mathcal F_0$ such that $f(x)=y$
then $z\in H_f$ and  if there is $f\in \mathcal F_1$ such that
$f(y)=x$ then $z\in {}_f H$.
\end{proof}

\subsubsection{Varying the Lipschitz constants}

\begin{lemma}\label{changingslope} Let $a,b>0$.  Then
$\cov(\Lip_{a,b})=\cov(\Lip_{2\cdot a,{\frac b 2}})$.
\end{lemma}

\begin{proof}
For $i\in 2$ let $X_i$ be the set of all sequences in $2^\omega$
starting with~$i$.  Let $h_i:2^\omega\to X_i$ be the homeomorphism
mapping $x$ to $(i,x(0),x(1),\dots)$.

Let $f:2^\omega\to 2^\omega$ be a Lipschitz function of constant
$a$. For $i\in 2$ let $f^i:2^\omega\to 2^\omega$ be a function
which is equal to $f\circ h_i^{-1}$ on $X_i$ and constant on
$X_{1-i}$ such that $f^i$ is Lipschitz of constant $2\cdot a$.
(For example, we can choose the constant value of $f^i$ on
$X_{1-i}$ to be $f((\overline{1-i}))$ where $(\overline{1-i})$
denotes the constant sequence with value~$1-i$.)

If $f:2^\omega\to 2^\omega$ is a Lipschitz function of constant
$b$, then for $i\in 2$ let $f_i:=h_i\circ f$.  $f_i$ is a
Lipschitz function of constant $\frac b 2$.

Now let $\mathcal F$ be a family of Lipschitz functions of
constant $a$ and $\mathcal G$ a family of Lipschitz functions of
constant~$b$. If $(2^\omega)^2=\bigcup\{f\cup
g^{-1}:f\in\mathcal F\wedge g\in\mathcal G\}$, then
$(2^\omega)^2=\bigcup\{f^i\cup g_i^{-1}:i\in 2\wedge
f\in\mathcal F\wedge g\in\mathcal G\}$.

It follows that $\cov(\Lip_{2\cdot a,{\frac b
2}})\leq\cov(\Lip_{a,b})$.  Now the lemma follows from the fact
that $\cov(\Lip_{a,b})$ is symmetric in $a$ and~$b$.
\end{proof}

\begin{lemma}\label{lipeqhm}
 Let $a,b>0$.  If there is $c\in\mathbb Z$ such that
$2^{c-1}\leq a$ and $2^{-c}\leq b$, then $\cov(\Lip_{a,b})=\hm$.
Otherwise $\cov(\Lip_{a,b})=2^{\aleph_0}$.
\end{lemma}

\begin{proof}
Let $a,b>0$ and assume that there is no $c\in\mathbb Z$ such that
$2^{c-1}\leq a$ and $2^{-c}\leq b$.  Let $c\in\mathbb Z$ be
maximal with $2^{c-1}\leq a$.  Then $2^{-c}>b$ and therefore
$b\cdot 2^c<1$.  $2^{c-1}\leq a$ is equivalent to $a\cdot
2^{-c}\geq {\frac 1 2}$, and since $c$ is maximal, we have $a\cdot
2^{-c}<1$.  By Lemma \ref{changingslope},
$\cov(\Lip_{a,b})=\cov(\Lip_{2^{-c}\cdot a,2^c\cdot b})$.  But
even the diagonal in $(2^\omega)^2$ cannot be covered
by less than $2^{\aleph_0}$ Lipschitz functions of constant~$<1$.

Now suppose there is $c\in\mathbb Z$ such that $2^{c-1}\leq a$ and
$2^c\leq b$.  By Lemma \ref{changingslope} we may assume $a\geq
1$ and $b\geq {\frac 1 2}$, hence $\cov(\Lip_{a,b})\ge
\cov(\Lip_{1,{\frac 1  2}})=\hm$.
\end{proof}

Let $\Lip$ be the $\sigma$-ideal on $2^\omega$ generated by
$\bigcup_{a>0}\Lip_{a,a}$, i.e., the $\sigma$-ideal generated by
all Lipschitz functions and their inverses.

\begin{thm}
$\hm=\cov(\Lip)$
\end{thm}

\begin{proof}
Clearly, $\cov(\Lip)\le \cov(\Lip_{1,{\frac 1 2}})$.  Thus, it
follows from Lemma \ref{limitedlipeqhm} that
$\cov(\Lip)\leq\hm$.

Now we prove the  converse inequality $\hm\leq\cov(\Lip)$.  We
define a coloring $c:[2^\omega\times 2^\omega]^2\to\mathcal P(2)$
as follows.

Let $0\in c((x_0,y_0),(x_1,y_1))$ iff there is a Lipschitz function of
constant $1$ containing both $(x_0,y_0)$ and $(x_1,y_1)$, i.e., if the
slope determined by $(x_0,y_0)$ and $(x_1,y_1)$ is $\leq 1$ or
equivalently, if $x_0$ and $x_1$ do not split after $y_0$ and~$y_1$.

Let $1\in c((x_0,y_0),(x_1,y_1))$ iff there is a Lipschitz function of
constant $1$ containing both $(y_0,x_0)$ and $(y_1,x_1)$, i.e., if
$y_0$ and $y_1$ do not split after $x_0$ and~$x_1$.

It is clear that $c$ is continuous and the color $\emptyset$ does not
occur.  We construct a (nonempty) perfect set $X\subseteq
(2^\omega)^2$ with the following properties:

\begin{itemize}\item[(i)] $c\restriction X$ only takes the values $\{0\}$
and $\{1\}$.
\item[(ii)] $c\restriction X$ is reduced.
\item[(iii)] For every Lipschitz function $f:2^\omega\to 2^\omega$,
$f\cap X$ and $f^{-1}\cap X$ are the unions of finitely many
$c$-homogeneous sets.
\end{itemize}

If we can construct~$X$, we are done.  This is because by (iii), every
family $\mathcal F$ of Lipschitz functions that
covers $(2^\omega)^2$ induces a family $\mathcal H$ of size at most
$\card{\mathcal F}$ that covers $X$ and consists of $c$-homogeneous
sets.  By (i) and (ii), we have $\hm\leq\card{\mathcal H}$ and thus
$\hm\leq\card{\mathcal F}$.

The required $X$ will be chosen to be (the graph of) a
homeomorphism between two perfect subsets of $2^\omega$.  For its
construction, partition $\omega$ into countably many intervals
$I_i$, $i\in\omega$, such that the length of every $I_i$ is at
least $i$ and the elements of $I_i$ are below the elements of
$I_j$ for~$i<j$.  For every $i\in\omega$ let $n_i$ denote the
first element of~$I_i$.

Let $T_0$ be a perfect subtree of $2^{<\omega}$ that fully splits at
all the levels of height $n_i$ for even $i$ and does not split at any
other level.  Let $T_1$ be a perfect subtree of $2^{<\omega}$ that
fully splits at every level of height $n_i$ for odd $i$ and does not
split anywhere else.

Let $X$ be the (graph of the) natural (order preserving) homeomorphism
between $[T_0]$ and $[T_1]$.  Clearly $X$ is closed and satisfies (i)
and (ii).  It remains to show (iii).

Let $f:2^\omega\to 2^\omega$ be a Lipschitz function.  Choose
$i\in\omega$ so that the Lipschitz constant of $f$ is below $2^{n_i}$.
$T_0$ is the union of finitely many perfect subtrees
$T_0^1,\dots,T_0^m$ that have no splittings below level~$n_i$.  For
$k\in\{1,\dots,m\}$ let $X_k:=X\cap([T_0^k]\times 2^\omega)$.  It is
straightforward to check that for all $k\in\{1,\dots,m\}$, $f\cap X_k$
is $c$-homogeneous of color $\{0\}$.

Similarly, the intersection of every inverse of a Lipschitz function
with $X$ is the union of finitely many $c$-homogeneous sets of color
$\{1\}$.  This shows (iii) and therefore finishes the proof of the
theorem.
\end{proof}

\subsubsection{Covering $(\om^\om)^2$ and $\R^2$ by Lipschitz functions}
  We
generalize our notation $\Lip_{a,b}$ to metric spaces~$X$.  For a
metric space $X$ let $\Lip_{a,b}(X)$ be the $\sigma$-ideal on
$X\times X$ generated by the Lipschitz functions of constant $a$
and the reflections of Lipschitz functions of constant~$b$.
$\Lip(X)$ denotes the $\sigma$-ideal generated by the union of
all the ideals $\Lip_{a,b}(X)$.

Recall that $\hm(c_{\parity})=\hm$.  It is easily checked that the
main arguments for the correspondence between Lipschitz functions on
$2^\omega$ and $c_{\min}$-homogeneous sets also go through for
$\omega^\omega$ and $c_{\parity}$.  This shows

\begin{corollary}\label{bairelip} For $X=2^\omega$ and $X=\omega^\omega$
we have
$$\cov(\Lip(X))=\cov(\Lip_{1,\frac 1 2}(X))=\hm.$$
\end{corollary}

At the very moment we do not know the exact relation between the
cardinal invariants mentioned above and $\cov(\Lip(\mathbb R))$.
However, we can say something:

\begin{rem} $\hm\leq\cov(\Lip(\mathbb R))$
\end{rem}

\begin{proof} The argument is similar to the argument in the proof of
Theorem \ref{lipeqhm}.

It is not difficult to construct a topological embedding
$e:2^\omega\to\mathbb R^2$ such that for any two distinct points
$x,y\in 2^\omega$ the slope determined by $e(x)$ and $e(y)$ is positive
and $\ge\Delta(x,y)$ if $\Delta(x,y)$ is even and $\le{\frac
{1}{\Delta(x,y)}}$
if $\Delta(x,y)$ is odd.

If $f:\mathbb R\to\mathbb R$ is Lipschitz, then $e^{-1}[f]$
is a finite union of $c_{\min}$-homogeneous sets of color $1$
and $e^{-1}[f^{-1}]$ is a finite union of $c_{\min}$-homogeneous sets of
color~$0$.  A covering family of Lipschitz real functions induces a
covering family of no greater size of $c_{\min}$-homogeneous subsets of
$2^\om$.  This
implies $\hm\leq\cov(\Lip(\mathbb R))$.
\end{proof}

\subsection{Covering squares by continuous functions}

After having established the equality
$\hm=\cov(\Lip_{1,{\frac 1 2}})=\cov(\Lip)$ and the fact that the Lipschitz
constants can be varied to some extent without changing
$\cov(\Lip_{a,b})$, it is natural to ask what happens if we replace
the Lipschitz functions by continuous functions.

For a topological space $X$ let $\Cont(X)$ denote the
$\sigma$-ideal on $X\times X$ generated by the continuous
functions from $X$ to $X$ and their inverses.  $\Cont$ is
$\Cont(2^\omega)$.  Obviously, $\Lip_{a,b}\subseteq\Cont$ for all
$a,b>0$. Theorem \ref{coveringbyfunctions} implies that
$\cov(\Cont)^+\geq 2^{\aleph_0}$.  The same is of course true for
$\hm$.  The question is whether $\cov(\Cont)$ can be smaller than
$\hm$. This will be answered in the next Section \ref{lipcontsep}.

Very often cardinal invariants of $\sigma$-ideals on Polish
spaces do not depend on the particular space the ideal is defined
on. This is not true for $\cov(\Cont(X))$. While
$\cov(\Cont(2^\omega))$ is consistently smaller than
$2^{\aleph_0}$,  the fact that every continuous function
from a connected space to a zero-dimensional space is constant implies
easily that if $X$ is the disjoint union of
$\mathbb R$ and $2^\omega$, then $\cov(\Cont(X))=2^{\aleph_0}$.

\subsubsection{The crucial inequality}
We show that $\cov(\cont(X))$ is  the same for $X=2^\omega$,
$X=\omega^\omega$, and $X=\mathbb R$. The proof of this fact
depends on the following perhaps surprising Theorem. The proof
below is the only proof in Part I which uses mathematical logic
techniques.

\begin{thm}\label{geqd} $\cov(\Cont(2^\omega))\geq\mathfrak d$,
where $\mathfrak d$ is the dominating number.
\end{thm}

\begin{proof} Let $\mathcal F$ be a family of continuous functions from
$2^\omega$ to $2^\omega$. Let $M$ be an elementary submodel of a
sufficiently large initial segment of the universe with Skolem
functions such that $\mathcal F\subseteq M$ and
$\card{M}=\card{\mathcal F}$.

Suppose $\card{\mathcal F}<\mathfrak d$. Then there is a function
$x\in 2^\omega\setminus M$. Let $M[x]$ denote the Skolem hull of
$M\cup\{x\}$. Since $\card{M[x]}=\card{\mathcal F}<\mathfrak d$,
there is $y:\omega\to\omega$ such that $y$ is not eventually
dominated by any function in $\omega^\omega\cap M[x]$.

Let $g:\omega^\omega\to 2^\omega$ be the natural embedding, i.e., the
one induced by the mapping that maps $n\in\omega$ to the sequence of
zeros of length $n$ followed by a single one.  Clearly, $g\in M$.

No $f\in\mathcal F$ maps $x$ to~$g(y)$, since such an $f$ would be an
element of $M[x]$ and therefore $y=g^{-1}(f(x))$ would be an element
of~$M[x]$.

Assume that there is $f\in\mathcal F$ such that $f(g(y))=x$.  Let
$h:=f\circ g$.  Then $h\in M$.

We work in $M[x]$ for a moment.  Since no function in $M[x]$
eventually dominates all functions from $D:=h^{-1}(x)$, by
elementarity, there is no
function at all which eventually dominates every function in~$D$.  In
other words, $D$ is unbounded.

A result of Kechris \cite{kechris} says
that every unbounded and closed set $D\subseteq\omega^\omega$
satisfies $D=A\cup P$, $A\cap P=\emptyset$ where $A$ is bounded, i.e.,
a single function eventually dominates all functions in~$A$, and
$P$ is superperfect, i.e., for all $s\in T(P)$ there is $t\in T(P)$
such that $s\subseteq t$ and $\succ_{T(P)}(t)$ is infinite. Since
$M[x]$ is elementary and $D=h^{-1}(x)\in M[x]$ is unbounded and
closed, there exist $A,P$ as above in~$M[x]$.

Now consider the set $B$ of all branches of $T(P)$ that do not meet
any node with infinitely many immediate successors in~$T(P)$.  It is
easy to see that $B$ is compact and thus bounded. Since $B$ is
definable in $M[x]$ and thus bounded by a function in~$M[x]$, $y$
cannot be an element of~$B$.  It follows that $y$ has an initial
segment $s\in T(P)$ such that $\succ_{T(P)}(s)$ is infinite.

As before, for $t\in 2^{<\omega}$ let $O_t$ denote the basic open
subset of $2^\omega$ consisting of all extensions of~$t$.

\begin{claim}
There is $t\in 2^{<\omega}$ such that for all $i\in 2$,
$T(h^{-1}[O_{t\frow i}])\cap\succ_{\omega^{<\omega}}(s)$ is infinite.
\end{claim}

\begin{proof}
Suppose not.  Then for all $t\in 2^{<\omega}$ at most one of the sets
$T(h^{-1}[O_{t\frow i}])\cap\succ_{\omega^{<\omega}}(s)$, $i\in 2$, is
infinite.  It follows that there is a unique $z\in 2^\omega$ such that
for all $n\in\omega$, $T(h^{-1}[O_{z\restriction
n}])\cap\succ_{\omega^{<\omega}}(s)$ is infinite.  Since $h\in M$,
$z\in M$.

But $x$ satisfies the definition of~$z$.  Thus $x=z$ and $x\in M$.  A
contradiction.
\end{proof}

Now let $t$ be as guaranteed by the claim.  Since $f$ is uniformly
continuous, there is $n\in\omega$ such that: \\
 for all $i\in 2$ and all
$r\in 2^n\cap T(f^{-1}[O_{t\frow i}])$ we have:   $f[O_r]\subseteq
O_{t\frow i}$.

For sufficiently large $m\in\omega$ we have for all
$a,b\in\omega^\omega$,
$$s\subseteq a\cap
b\wedge
\min\{a(\dom(s)),b(\dom(s))\}>m \;\Rightarrow \; g(a)\restriction
n=g(b)\restriction n,$$
from which it follows that for sufficiently large
$m\in\omega$:
\begin{multline*}
s\subseteq a\cap
b\wedge
\min\{a(\dom(s)),b(\dom(s))\}>m \;\Rightarrow\\ t\subseteq
h(a)\wedge
\bigl(h(a)\restriction(\dom(t)+1)=h(b)\restriction(\dom(t)+1)\bigr).
\end{multline*}

But this contradicts the choice of $t$ and hence no $f\in\mathcal F$
maps $g(y)$ to~$x$.  Thus, $(2^\omega)^2$ is not covered by $\mathcal F$.
This shows $\cov(\Cont(2^\omega))\geq\mathfrak d$.
\end{proof}

It should be pointed out that the curious use of two new reals over
$M$ in the proof of Theorem \ref{geqd} is really necessary. It can be
shown that after adding a Miller real, which is unbounded, $2^\omega$
is covered by the $c_{\min}$-homogeneous sets coded in the ground
model.  In particular, after adding one Miller real, $(2^\omega)^2$
is covered by the continuous functions coded in the ground
model.  The proof of Theorem \ref{geqd} shows that
this is not the case after adding two Miller reals or even any new
real and then a Miller real over it.

\relax From Theorem \ref{geqd} we get

\begin{thm}\label{ZFCequalities}
$\cov(\Cont(\omega^\omega))=\cov(\Cont(\mathbb
R))=\cov(\Cont(2^\omega))$
\end{thm}

The proof of this theorem uses the following easy observation.

\begin{lemma}\label{context} Let $C$ be a compact subset of
$\omega^\omega$ and let $f:C\to\omega^\omega$ be continuous. Then $f$
can be continuously extended to all of $\omega^\omega$.
\end{lemma}

\begin{proof} Consider $C$ as a subset of $(\omega+1)^\omega$.  The latter
space is homeomorphic to $2^\omega$.  $f[C]$ is bounded and therefore
in $\omega^\omega$ there is a copy of $2^\omega$ including~$f[C]$. The
lemma now follows from the well-known fact that every continuous
mapping from a closed subset of a Boolean space to $2^\omega$ can be
continuously extended to the whole space (which follows from
$2^\omega$ being the Stone space of a free Boolean algebra).
\end{proof}

\begin{proof}[Proof of Theorem \ref{ZFCequalities}]
We first show
$\cov(\Cont(2^\omega))\leq\cov(\Cont(\omega^\omega))$. Let
$f:\omega^\omega\to\omega^\omega$ be continuous.  Then
$f^{-1}[2^\omega]$ is closed and thus $A:=f^{-1}[2^\omega]\cap
2^\omega$ is a closed subset of $2^\omega$.  We can now extend
$f\restriction A$ to a continuous function $\overline f:2^\omega\to
2^\omega$ by the same argument as in the proof of Lemma \ref{context}.

This shows that a family $\mathcal F\subseteq\Cont(\omega^\omega))$
covering $(\omega^\omega)^2$ gives rise to a covering
family of no greater size in $\Cont(2^\omega)$ and thus,
$\cov(\Cont(2^\omega))\leq\cov(\Cont(\omega^\omega))$.

The same argument goes through for $\mathbb R$ instead of
$\omega^\omega$, using the Tietze-Urysohn theorem.

Now observe that $\omega^\omega$ can be covered by $\mathfrak d$ copies
of $2^\omega$ since $\mathfrak d$ is the covering number of the ideal
of bounded subset of $\omega^\omega$.  Let $\mathcal C$ be a
collection of size $\mathfrak d$ of copies of $2^\omega$ covering
$\omega^\omega$.

To each pair $(C,D)\in{\mathcal C}\times{\mathcal C}$
assign a family $\mathcal F_{C,D}$ of size
$\cov(\Cont(2^\omega))$ of continuous functions on $\omega^\omega$
such that $C\times D\subseteq\bigcup\{f\cup f^{-1}:f\in\mathcal
F_{C,D}\}$. This is possible by Lemma \ref{context}. Let $\mathcal
F:=\bigcup\{\mathcal F_{C,D}:C,D\in\mathcal C\}$. Now
$(\omega^\omega)^2=\bigcup\{f\cup f^{-1}:f\in\mathcal F\}$ and the
size of $\mathcal F$ is $\max(\cov(\Cont(2^\omega)),\mathfrak
d)=\cov(\Cont(2^\omega))$.  The last equality is Theorem \ref{geqd}.

Again, same argument works for $\cov(\Cont(\mathbb R))$ as well since
$\mathbb R$ is just $\omega^\omega$ (the irrationals) together with
countably many additional points (the rationals) and therefore also
can be covered by $\mathfrak d$ copies of $2^\omega$. We again use
the Tietze-Urysohn theorem to extend continuous mappings defined on
closed subspaces of $\mathbb R$.
\end{proof}

\bigbreak
\bigbreak
\begin{center}
\Large{\sc Part II: Independence results}
\end{center}

In the second part of the paper we show that any two rows in
Diagrmam~1 can be
separated. We shall prove that every assignment of $\aleph_1$-s
and $\aleph_2$-s to the diagram which is consistent with the
arrows is realized in a model of set theory.

We provide two new forcing notions. One for separating
$\hm(c_{\min})$ from $\hm(c_{\max})$ and the other for separating
$\cov(\Cont(2^\omega))$ from $\cov(\Lip(2^\omega))$. We force
over models of CH with countable support iterations of Axiom A
forcing  notions (see \cite{axioma}) of size $\aleph_1$ which add new
reals.  Thus, no cardinals are collapsed and in the resulting
models the continuum is $\aleph_2$.

Theorem \ref{coveringbyfunctions} implies that if the continuum
is a limit cardinal,  all three numbers above are equal to the
continuum.  In fact, it is very easy to make
$\cov(\Cont(2^\omega))$ equal to the continuum.

Let $M$ be a model of set theory and assume that $\mathcal F\in
M$ is a family of continuous functions on $2^\omega$.  If $x,y\in
2^\omega$ are generic over $M$ and independent in the sense that
$x\not\in M[y]$ and $y\not\in M[x]$, then no $f\in\mathcal F$ can
cover $(x,y)$. It follows that after forcing with a large
product of some sort in order to increase the continuum one ends
up with a model of set theory where $\cov(\Cont(2^\omega))$ is
the continuum.  In particular, after forcing with the measure
algebra over $2^{\aleph_2}$ over a model of CH, one obtains a
model (the Solovay model) in which $\d=\aleph_1$ (since the
ground model elements of $\omega^\omega$ dominate the new
elements) and $\cov(\Cont(2^\omega))=\aleph_2$.

In \cite{GKKS} it was shown that in the Sacks model all
homogeneity numbers of reduced continuous pair-colorings on
Polish spaces are equal to $\aleph_1<2^{\aleph_0}$.  It follows
that $\hm(c_{\max})$, $\hm(c_{\min})$, and $\cov(\Cont(\mathbb
R))$ are small in the Sacks model.

There is a  natural forcing $\mathbb P_{c_{\min}}$ for  separating
$\hm(c_{\min})$ from the numbers below it: forcing with Borel
subsets of $2^{\om}$ which are positive with respect to the
$\sigma$-ideal $J_{\min}$ generated over $2^\om$ by all
$c_{\min}$-homogeneous sets.  We show that the countable support
iteration of this forcing of length $\omega_2$ produces a model of
$\hm(c_{\min})=\aleph_2$ and $\cov(\Cont(2^\omega))=\aleph_1$. In
this model it  holds that covering $\R^2$ by Lipschitz functions
is strictly more difficult than covering $\R^2$ by continuous
functions.

 By a new (and yet unpublished)
theorem  of Zapletal, the existence of large cardinals
implies that $\mathbb P_{c_{\min}}$ is optimal for enlarging
$\hm(c_{\min})$ in the sense that it does not enlarge numbers
which are consistently smaller than $\hm(c_{\min})$. Assuming
large cardinals, Shelah and Zapletal proved recently that for
every reasonably defined $\sigma$-ideal $I$ on the reals whose
covering number is provably $\ge\hm(c_{\min})$,
 the uniformity of
$I$  (i.e., the smallest size of a set not in~$I$) is at most
$\aleph_3$.  (The uniformity of $J_{\min}$ is at most $\aleph_2$.)

The analogous natural forcing for increasing $\hm(c_{\max})$ is,
however, not only not optimal, but actually \emph{increases} the
smaller $\hm(c_{\min})$. So another forcing has to be used for
separating $\hm(c_{\min})$ from $\hm(c_{\max})$.

We design a new tree-forcing notion $\mathbb P_{c_{\max}}$  for
increasing $\hm(c_{\max})$ while leaving $\hm(c_{\min})$ small.
The tree-combinatorics required for this forcing stems from   a
new result of Noga Alon about a Ramsey connection between perfect
graphs and random graphs \cite{Alon} (which Alon proved for this
purpose). The countable support iteration of length $\om_2$ of
$P_{c_{\max}}$  produces a model of set theory in which
$\hm(c_{\min})=\aleph_1$ and $\hm(c_{\max})=\aleph_2$.

\section{Consistency of  $\hm(c_{\min})<\hm(c_{\max})$}

\subsection{The $c_{\max}$-forcing} We are looking for a notion of
forcing which adds a real that avoids all the $c_{\max}$-homogeneous
sets in the ground model but does not increase $\hm$ when iterated.

\begin{defn}  For a pair-coloring $c$ of a finite set with two
colors let $\norm(c)$ denote the greatest $n\in\omega$ for which
$\chi_{\random}\restriction n\leq c$.

For a subtree $p\subseteq T(\prod_{n\in\omega}(n+1))$ and $t\in p$
let
$c_{t,p}:=\chi_{\random}\restriction\{i\in\omega:t\frow
i\in\succ_p(t)\}$.
(See \ref{notation} for notation.)   % mg May 11

Let
\begin{align*}\mathbb P_{c_{\max}}:=\Bigl\{ p \su
T\bigl(\prod_{n\in\omega}(n+1)\bigr):\quad &p \text{ is a tree and
}\\
&\forall t\in p\forall n\in\omega\exists s\in p
(s\supseteq t\wedge\norm(c_{s,p})\geq n)\Bigr\}
\end{align*}

The order on $\mathbb P_{c_{\max}}$ is set-inclusion.
\end{defn}

In the following we write just $\mathbb P$ for $\mathbb
P_{c_{\max}}$.  For a condition $p\in \mathbb P$ and $t\in p$,
let
$p_t=\{s\in p: s \subseteq t \vee t\subseteq s\}$, and call $p_t$ the
condition $p$ \emph{below}~$t$.
It is clear that
$p_t\in \mathbb P$  for $p\in \mathbb P$ and $t\in p$. If
$G\subseteq
\mathbb P$ is a generic filter over a ground model~$M$, the
\emph{generic real added by $\mathbb P$} is
the unique element of      % mg May 11
$\bigcap \{[p]:p\in G\}$.

\begin{claim}\label{gennotina}  The generic real added by
$\mathbb P$ avoids all $c_{\max}$-homogeneous sets in the ground
model.
\end{claim}

\begin{proof} Suppose  that $A\su 2^\omega$ is $c_{\max}$-homogeneous, say
with color~$0$, and $A\in M$.
Let $p\in
\mathbb P$ be arbitrary. Choose
$s\in p$  with $t,t'\in \Succ_p(s)$ satisfying
$\random(t,t')=1$. Since at least one of $[p_t]$, $[p_{t'}]$ has
empty intersection with~$A$, assume without loss of generality that
$[p_t]\cap A=\emptyset$. Now
$p_t\le p$ is a condition in $\mathbb P$ which forces that the generic
real is not in~$A$. Thus, the set of conditions forcing that the
generic real is not in $A$ is dense and belongs to~$M$, hence the
generic real is not in~$A$.
\end{proof}

\begin{lemma}\label{onereal} Let $G$ be $\mathbb P$-generic over the
ground model~$M$.  Then for each $x\in(2^\omega)^{M[G]}$ there is a
tree $T\in M$ such that $[T]$ is $\parity$-homogeneous and $x\in[T]$.
\end{lemma}

For the proof of this lemma we use the following result of Noga Alon
\cite{Alon} that was proved especially for this purpose.

\begin{lemma}\label{binarycoloring} Let $n\in\omega$ and $c:[n]^2\to 2$.
Then there is $N\in\omega$ and $C:[N]^2\to 2$ such that whenever
$e:N\to 2^\omega$ is 1-1, then there is a $c_{\min}$-homogeneous set
$H\subseteq 2^\omega$ such that $c\leq C\restriction e^{-1}[H]$.
\end{lemma}

\begin{proof}[Proof of Lemma \ref{onereal}]
Let $\dot{x}$ be a name for a new element of $2^\omega$ and let
$p\in\mathbb P$.  Since $\dot{x}$ is a name for a new real, we may
assume, by passing to stronger condition if necessary, that for each
splitting node $s\in p$ and all $t,t^\prime\in\succ_p(s)$ with
$t\not=t^\prime$, the initial segments of $\dot{x}$ decided by $p_t$
and $p_{t^\prime}$ are incompatible.

We may assume that for some $k_s\in\omega$, each~$p_t$,
$t\in\succ_p(s)$, decides an initial segment of $\dot{x}$ of length
$k_s$ and that the decisions of the $p_t$'s on $\dot{x}$ are
pairwise incompatible when restricted to~$k_s$.  In other words, for
each splitting node $s$ of $p$ we have an embedding $e_s:\succ_p(s)\to
2^{k_s}$
with $p_t \forces \dot x \supseteq e_s(t)$.

Now Lemma \ref{binarycoloring} implies that we can thin out $p$ to a
condition $q$ such that for each splitting node $s$ of~$q$,
$e_s[\succ_q(s)]$ is a $c_{\min}$-homogeneous subset of $2^{k_s}$ of
some color $i_s\in 2$.

Thinning out $q$ further if necessary, we may assume that
\begin{itemize}\item[$(*)$]  whenever $s$ and
$t$ are splitting nodes of $q$ and $s\subsetneqq t$, then
$\norm(c_{s,q})<\norm(c_{t,q})$.
\end{itemize}

Now either $q$ has a cofinal set of splitting nodes $s$ with $i_s=0$,
or there is a node $s\in q$ such that for all splitting nodes $t\in q$
with $s\subseteq t$, $i_t=1$.  In the first case, we can thin out $q$
to a condition $r$ such that for all splitting nodes $s$ of~$r$,
$i_s=0$.  The property~$(*)$
makes sure that $r$ will be a condition. In the second
case we can put $r:=q_s$ and get a condition such that for all
splitting nodes $s$ we have $i_s=1$.

Finally let $T_r:=\{s\in 2^{<\omega}:\exists r^\prime\leq r
(r^\prime\forces s\subseteq\dot{x})\}$ be the {\em tree of
$r$-possibilities for $\dot{x}$}.  Clearly $r$ forces $\dot{x}$ to be
a branch of~$T_r$. By the construction of~$r$, $[T_r]$ is
$c_{\min}$-homogeneous.
\end{proof}

\subsection{Iteration}

In this section we show that after forcing with a countable support
iteration of the $c_{\max}$-forcing, all the new reals ($\in
2^\omega$) are covered by $c_{\min}$-homogeneous sets in the ground
model. This implies that after forcing with a countable support
iteration of $\mathbb P$ of length $\omega_2$ over a model of CH, we
obtain a model of set theory in which $\hm=\aleph_1$ but
$\hm(c_{\max})=\aleph_2$.  The latter statement follows from
Claim \ref{gennotina}.

\subsubsection{A preliminary lemma} Our strategy is the following:
For an ordinal $\alpha$ let $\mathbb P_\alpha$ denote the
countable support iteration of $\mathbb P$ of length $\alpha$.
Let $\dot{x}$ an $\mathbb P_{\omega_2}$-name for a new element of
$2^\omega$.  We may assume that there is $\alpha<\omega_2$ such
that $\dot{x}$ is an $\mathbb P_\alpha$-name for a real not added
at any proper initial stage of the iteration $\mathbb P_\alpha$.
Let $q$ be a condition in $\mathbb P_{\alpha}$. Recall the
definition of $T_q\subseteq 2^{<\omega}$ from the proof of
Lemma \ref{onereal}:
$$T_q=\{s\in 2^{<\omega}:\exists q^\prime\leq q(q^\prime\forces
s\subseteq\dot{x})\}.$$

For each $p\in\mathbb P_{\alpha}$ we will construct a condition $q\leq
p$ such that $[T_q]$ is $c_{\min}$-homogeneous.  The next lemma tells
us how to choose the color of $[T_q]$ if $\dot{x}$ is added in a limit
step.  That is, we can decrease $p$ such that it becomes an element of
one of the sets~$E_i$, $i\in 2$, defined below. If $p\in E_i$, we can
build $q\leq p$ such that $[T_q]$ is $c_{\min}$-homogeneous of color
$i$.

Let us fix some notation.  If $\mathbb Q$ is any forcing notion and
$\dot{y}$ is a $\mathbb Q$-name for a new element of $2^\omega$ let
$y[p]$ be the maximal element of $2^{<\omega}$ such that $p\forces
y[p]\subseteq\dot{y}$.  $y[p]$ exists since $\dot{y}$ is a name for a
new real.

For $i\in 2$ let
\begin{multline*}
E_i:=\{p\in\mathbb P_\alpha:\forall\beta<\alpha\forall q\leq p\exists
q^\prime\leq q\exists q_0,q_1\in \mathbb P_{\beta, \alpha}
\\ (q^\prime\restriction\beta\forces
q_0,q_1\leq
q^\prime\restriction[\beta,\alpha)\wedge\parity(x[q_0],x[q_1])=i)\}.
\end{multline*}

Recall that $\parity(s,t)\in 2$ implies that $s$ and $t$ are
incompatible, i.e., $s\perp t$.

\begin{lemma}\label{densesets} $E_0$ and $E_1$
are open and $E_0\cup E_1$ is dense in $\mathbb P_\alpha$.
\end{lemma}

This lemma is true for all forcing iterations, not only of variations
of Sacks forcing. We do not even use the countable supports.

\begin{proof}[Proof of Lemma \ref{densesets}] Let us start with
\begin{claim}\label{changingordinals}
Let $\beta<\alpha$ and let $q\in\mathbb P_\alpha$ be such that for
some $i\in 2$ there are $q_0$ and $q_1$ such that
$$q\restriction\beta\forces q_0,q_1\leq
q\restriction[\beta,\alpha)\wedge\parity(x[q_0],x[q_1])=i.$$ Let
$\gamma<\beta$.  Then there are $q^\prime\leq q$ and $q^\prime_0$ and
$q^\prime_1$ such that $$q^\prime\restriction\gamma\forces
q^\prime_0,q^\prime_1\leq
q^\prime\restriction[\gamma,\alpha)\wedge\parity(x[q^\prime_0],x[q^\prime_1])=i.$$
\end{claim}

To see this, let $q^\prime\leq q$ be such that
$q^\prime\restriction[\beta,\alpha)=q\restriction[\beta,\alpha)$ and
$q^\prime\restriction\beta$ decides $x[q_0]$ and $x[q_1]$.  For $j\in
2$ let $q_j^\prime:=(q^\prime\restriction[\gamma,\beta))\frow
q_j$. Now
$q^\prime$, $q^\prime_0$, and $q^\prime_1$ work for the claim.

For the proof of the lemma let $p\in\mathbb P_\alpha$. Suppose
$p\not\in E_0$.  We show that $p$ has an extension in~$E_1$.  Since
$p\not\in E_0$, there are $\gamma<\alpha$ and $q\leq p$ such that for
all $q^\prime\leq q$ and any two sequences $q_0$ and $q_1$ for names
of conditions, if $q^\prime\restriction\gamma\forces q_0,q_1\leq
q^\prime\restriction[\gamma,\alpha)$, then
$q^\prime\restriction\gamma\not\forces\parity(x[q_0],x[q_1])=0$.  We
are done if we can show

\begin{claim} $q\in E_1$. \end{claim}

Let $r\leq q$ and $\beta<\alpha$.  Note that by Claim
\ref{changingordinals}, the sets $E_i$ are not changed if in the
definition we replace ``$\forall\beta<\alpha$'' by ``for cofinally
many $\beta<\alpha$''.  Thus we may assume $\beta\geq\gamma$.

Since we assumed that $\dot{x}$ is not added in a proper initial stage
of the iteration (before $\alpha$), there are $q_0$ and $q_1$ such
that
$$r\restriction\beta\forces q_0,q_1\leq
r\restriction[\beta,\alpha)\wedge x[q_0]\perp x[q_1].$$ Decreasing
$r\restriction\beta$ if necessary, we may assume that
$r\restriction\beta$ decides $\parity(x[q_0],x[q_1])$ to be $i\in 2$.

By Claim \ref{changingordinals}, there are $r^\prime\leq r$ and $r_0$
and $r_1$ such that
$$r^\prime\restriction\gamma\forces r_0,r_1\leq
r^\prime\restriction[\gamma,\alpha)\wedge\parity(x[r_0],x[r_1])=i.$$
By the choice of~$q$, $i\not=0$.  Thus~$i=1$. This shows $q\in E_1$.
\end{proof}

\subsubsection{Some forcing notation}\label{forcingnotation}

For $n\in\omega$ and $p\in\mathbb P$ let $p^n$ be the set of all minimal
$t\in p$ such that $\norm(c_{t,p})>n$. For $p,q\in\mathbb P$ we write
$q\leq_n p$ if $q\leq p$ and $p^n=q^n$.

A sequence $(p_n)_{n\in\omega}$ in $\mathbb P$ is a {\em fusion
sequence} if there is a nondecreasing unbounded function
$f:\omega\to\omega$ such that for all $n\in\omega$,
$p_{n+1}\leq_{f(n)}p_n$.  If $(p_n)_{n\in\omega}$ is a fusion
sequence, then $p_\omega=\bigcap_{n\in\omega}p_n$ is a condition in
$\mathbb P$, the {\em fusion} of the sequence.  In this definition,
the function $f$ is only added for technical convenience.  If we only
talk about the identity function instead of arbitrary~$f$, we arrive
at an essentially equivalent notion.

The idea behind fusion is that in $\mathbb P$, even though it is not
countably closed, lower bounds exist for suitably chosen countable
sequences.  All we have to do while inductively thinning out a
condition, is to leave splitting nodes with more and more complicated
colorings on their successors untouched.  This is exactly what we did,
although less formally, in the proof of Lemma \ref{onereal}. The
method can be extended to countable support iterations.

Let $\alpha$ be an ordinal. For $F\in[\alpha]^{<\aleph_0}$,
$\eta:F\to\omega$, and $p,q\in\mathbb P_\alpha$ let $q\leq_{F,\eta}p$
                 % mg May 11
if $q\leq p$ and for all $\beta\in F$, $q\restriction\beta\forces
q(\beta)\leq_{\eta(\beta)} p(\beta)$.  Roughly speaking,
$q\leq_{F,\eta}p$ means that on each coordinate from~$F$, $q$ is
$\leq_n$-below $p$ where $n$ is given by~$\eta$.

A sequence $(p_n)_{n\in\omega}$ of conditions in $\mathbb P_\alpha$ is
a {\em fusion sequence} if there is an increasing sequence
$(F_n)_{n\in\omega}$ of finite subsets of $\alpha$ and a sequence
$(\eta_n)_{n\in\omega}$ such that for all $n\in\omega$,
$\eta_n:F_n\to\omega$, $p_{n+1}\leq_{F_n,\eta_n}p_n$, for all
$\gamma\in F_n$ we have $\eta_n(\gamma)\leq\eta_{n+1}(\gamma)$, and
for all $\gamma\in\supt(p_n)$ there is $m\in\omega$ such that
$\gamma\in F_m$ and $\eta_m(\gamma)\geq n$.

This notion is precisely what is needed in countable support
iterations to get suitable fusions.  It essentially means that once we
have touched (i.e., decreased) a coordinate of~$p_0$, we have to build
a fusion sequence in that coordinate.

If $(p_n)_{n\in\omega}$ is a fusion sequence in $\mathbb P_\alpha$,
its fusion $p_\omega$ is defined inductively.  Let $F_\omega:=\bigcup
F_n$.

Suppose $p_\omega(\gamma)$ has been defined for all $\gamma<\beta$ for
some $\beta<\alpha$.  If $\beta\not\in F_\omega$, let
$p_\omega(\beta)$ be a name for $1_\mathbb P$.  If $\beta\in
F_\omega$, then $p_\omega\restriction\beta$ forces
$(p_n(\beta))_{n\in\omega}$ to be
a fusion sequence in $\mathbb P$.  Let $p_\omega(\beta)$ be a name for
the fusion of the $p_n(\beta)$'s.

\subsubsection{Keeping $\hm$ small}\label{hmsmall} Let $\dot{x}$
and $\alpha$ be as before.  The way to build a condition $q$ for
which $T_q$ is $c_{\min}$-homogeneous is the following: $q$ will
be the fusion of a fusion sequence $(p_n)_{n\in\omega}$ with
witness $(F_n,\eta_n)_{n\in\omega}$.  For each~$n$,
$(p_n,F_n,\eta_n)$ will determine a finite initial segment $T_n$
of~$T_q$.  We have to make sure that $T_q$ is the union of the
$T_n$ and that the $T_n$ are good enough to guarantee the
$c_{\min}$-homogeneity of~$[T_q]$.  The latter will be ensured by
the $(F_n,\eta_n)$-faithfulness of each~$p_n$, which is defined
below.

First we introduce some tools that help us to carry out the necessary
fusion arguments.

We call a condition $p\in\mathbb P$ {\em normal} if for every $s\in p$
with $n:=\card{\succ_p(s)}>1$, $c_{s,p}$ is isomorphic to
$c_{\random}\restriction n$ and moreover, if $t\in p$ is a minimal
proper extension of $s$ with more than one successor in~$p$, then
$\card{\succ_p(t)}=\card{\succ_p(s)}+1$.
 Thus, $s\in p^n $ iff $\card{\succ_p(s)} = n+1$.
% mg

Let $I:=T(\prod_{i\in \omega} (i+1))=
          \bigcup\{\prod_{i<n}(i+1):n\in\omega\}$ and $I_n:=\{\rho\in
I:\dom(\rho)=n\}$.  If $p\in\mathbb P$ is a normal condition, then
each $\rho\in I_n$ determines an element $s_\rho$ of~$p^n$. Let
$p*\rho:=p_{s_\rho}=\{t\in p:s_\rho\subseteq t\vee t\subseteq
s_\rho\}$.

A condition $q\in\mathbb P_\alpha$ is normal if for all
$\beta<\alpha$, $q\restriction\beta$ forces that $q(\beta)$ is normal.
Suppose $q\in\mathbb P_\alpha$ is a normal condition.  For
$F\in[\alpha]^{<\aleph_0}$, $\eta:F\to\omega$,
$\sigma\in\prod_{\gamma\in F}I_{\eta(\gamma)}$, and $q\in\mathbb
P_\alpha$ let $q*\sigma$ be defined as follows:

For $\gamma\in F$ let $(q*\sigma)(\gamma)$ be a name for a condition
in $\mathbb P$ such that $\forces_{\mathbb
P_\gamma}(q*\sigma)(\gamma)=q(\gamma)*\sigma(\gamma)$. For
$\gamma\in\alpha\setminus F$ let $(q*\sigma)(\gamma):=q(\gamma)$.

Now $(q*\sigma)_{\sigma\in\prod_{\gamma\in F}I_{\eta(\gamma)}}$
 is a finite maximal antichain below~$q$.  Consider
the tree $T$ generated by $\{x[q*\sigma]:\sigma\in\prod_{\gamma\in
F}I_{\eta(\gamma)}\}$.  If $q^\prime\leq_{F,\eta}q$, then
$T_{q^\prime}$ is an end-extension of~$T$.

It is clear that the normal conditions in $\mathbb P$ form a dense
subset and the same is true for $\mathbb P_\alpha$.  Therefore, from
now on all the conditions we consider are assumed to be normal.  We
have to be careful at one point, however.  Suppose $p\in\mathbb P$ is
a normal condition and we have constructed some $q\leq_n p$.  $q$ is
not necessarily normal.  But it is easy to see that there is some
$q^\prime\leq_n q$ which is normal.  We call the process of passing
from $q$ to $q^\prime$ {\em normalization at $n$}.  Normalization at
$n$ will be applied automatically without being mentioned whenever we
construct some $q\leq_n p$.

\begin{defn}\label{satdef}
Let $i\in 2$ and $\dot{x}$ be fixed.  For $F$ and $\eta$ as before, a
condition $q\in\mathbb P_\alpha$ is {\em $(F,\eta)$-faithful} if for
all $\sigma,\tau\in\prod_{\gamma\in F}I_{\eta(\gamma)}$ with
$\sigma\not=\tau$, $\parity(x[q*\sigma],x[q*\tau])=i$.
\end{defn}

Now we are ready to formulate the lemma that will allow us to handle
the case where $\dot{x}$ is added at a limit step of the iteration.

\begin{lemma}\label{limitstep}
Let $\alpha$ be a limit ordinal and let $\dot{x}$ be a $\mathbb
P_\alpha$-name for an element of $2^\omega$ which is not added by an
initial stage of the iteration.  Let~$F$, $\eta$, and $i$ be as in
Definition \ref{satdef} and suppose that $q\in\mathbb P_\alpha$ is
$(F,\eta)$-faithful.

a) Let $\beta\in\alpha\setminus F$ and let $F^\prime:=F\cup\{\beta\}$
and $\eta^\prime:=\eta\cup\{(\beta,0)\}$.  Then $q$ is
$(F^\prime,\eta^\prime)$-faithful.

b) Suppose $q\in E_i$.  Let $\beta\in F$ and let
$\eta^\prime:=\bigl(\eta\restriction
(F\setminus\{\beta\})\bigr)\cup\{(\beta,\eta(\beta)+1)\}$.  Then there is
$r\leq_{F,\eta}q$ such that $r$ is $(F,\eta^\prime)$-faithful.
\end{lemma}

\begin{proof}
a) follows immediately from the definitions.

For b) let $\delta:=\max(F)+1$ and $n:=\eta(\beta)$.
\begin{claim}  There is a condition
$q^\prime\leq_{F,\eta}q$ such that for each $\sigma\in\prod_{\gamma\in
F}I_{\eta(\gamma)}$ there are sequences
$q_{\sigma,0},\dots,q_{\sigma,n}$ of names for conditions such that
for all $k\leq n$,
$$q^\prime*\sigma\restriction\delta\forces q_{\sigma,k}\leq
q\restriction[\delta,\alpha),$$ $q^\prime*\sigma\restriction\delta$
decides $x[q_{\sigma,k}]$, and for all $l\leq n$ with $k\not=l$,
$$q^\prime*\sigma\restriction\delta\forces
\parity(x[q_{\sigma,k}],x[q_{\sigma,l}])=i.$$
\end{claim}
 For the proof of the claim, let $\{\sigma_1,\dots,\sigma_m\}$ be an
enumeration of $\prod_{\gamma\in F}I_{\eta(\gamma)}$. We build a
$\leq_{F,\eta}$-decreasing sequence $(q_j)_{j\leq m}$ such that
$q_0:=q$ and $q^\prime:=q_m$ works for the claim. As we construct
$q_j$, we find suitable $q_{\sigma_j,k}$ for all~$k<n$.

Let $j\in\{1,\dots,m\}$ and assume that $q_{j-1}$ has already been
constructed. Since $q\in E_i$ and $E_i$ is open, there are
$q^\prime_j\leq q_{j-1}*\sigma_j$ and sequences $q_{\sigma_j,0}$ and
$q^{\prime}_{\sigma_j,1}$ of names of conditions such that
$$q^\prime_j\restriction\delta\forces
q_{\sigma_j,0},q^\prime_{\sigma_j,1}\leq
q\restriction[\delta,\alpha)\wedge\parity
(x[q_{\sigma_j,0}],x[q^\prime_{\sigma_j,1}])=i.$$ Iterating this
process by splitting $q^\prime_{\sigma_j,1}$ into $q_{\sigma_j,1}$ and
$q^\prime_{\sigma_j,2}$ and so on and decreasing $q^\prime_j$, we
finally obtain $q^\prime_j\leq q_{j-1}$ and $q_{\sigma_j,k}$, $k\leq
n$, such that for all $k\leq n$.
$$q^\prime_j\restriction\delta\forces q_{\sigma_j,k}\leq
q\restriction[\delta,\alpha)$$ and for all $l\leq n$ with $l\not= k$,
$$q^\prime_j\restriction\delta\forces\parity
(x[q_{\sigma_j,k}],x[q_{\sigma_j,l}])=i.$$

We may assume that $q^\prime_j\restriction\delta$ decides
$x[q_{\sigma_j,k}]$ for all $k\leq n$. Let $q_j\leq_{F,\eta}q_{j-1}$
be such that
$q_j*\sigma_j\restriction\delta=q^\prime_j\restriction\delta$ and
$q_j\restriction[\delta,\alpha)=q\restriction[\delta,\alpha)$.  This
finishes the construction, and it is easy to check that it works.

Continuing the proof of lemma~\ref{limitstep},
let $q_{\sigma,k}$ and $q^\prime$
be as in the claim. For $\rho\in I_{\eta(\delta)}$ let $r^{\rho\frow
0},\dots,r^{\rho\frow n}$ be sequences of names for conditions such
that for all $k\leq n$ and all $\sigma\in\prod_{\gamma\in
F}I_{\eta(\gamma)}$ with $\sigma(\beta)=\rho$,
$$q^\prime*\sigma\restriction\delta\forces r^{\rho\frow
k}=q_{\sigma,k}.$$

Let $r$ be a sequence of names for conditions such that
$r\restriction\delta=q^\prime\restriction\delta$ and for all
$\sigma\in\prod_{\gamma\in F}I_{\eta^\prime(\gamma)}$,
$$q^\prime*\sigma\restriction\delta\forces
r\restriction[\delta,\alpha)= r^{\sigma(\beta)}.$$

Note that $r\leq_{F,\eta}q^\prime$ and thus $r\leq_{F,\eta}q$. It
follows from the construction that $r$ is $(F,\eta^\prime)$-faithful.
\end{proof}

A similar lemma is true if the new real is added in a successor step.

\begin{lemma}\label{successorstep}
Let $\alpha$ be a successor ordinal, say $\alpha=\delta+1$ and let
$\dot{x}$ be a $\mathbb P_\alpha$-name for an element of $2^\omega$
which is not added by an initial stage of the iteration.  Let~$F$,
$\eta$, and $i$ be as in Definition \ref{satdef} and suppose that
$q\in\mathbb P_\alpha$ is $(F,\eta)$-faithful.

a) Let $\beta\in\alpha\setminus F$ and let $F^\prime:=F\cup\{\beta\}$
and $\eta^\prime:=\eta\cup\{(\beta,0)\}$.  Then $q$ is
$(F^\prime,\eta^\prime)$-faithful.

b) Suppose $$q\restriction\delta\forces\mbox{``\/$[T_{q(\delta)}]$ is
$c_{\min}$-homogeneous of color $i$''}.$$ Let $\beta\in F$ and let
$\eta^\prime:=\eta\restriction
F\setminus\{\beta\}\cup\{(\beta,\eta(\beta)+1)\}$.  Then there is
$r\leq_{F,\eta}q$ such that $r$ is $(F,\eta^\prime)$-faithful.
\end{lemma}

\begin{proof} As in Lemma \ref{limitstep}, a) follows directly from the
definitions.

For the proof of b) we have to consider two cases.  First suppose
$\beta=\delta$.  In this case let $q^\prime$ be a name for a condition
in $\mathbb P$ such that for all $\sigma\in\prod_{\gamma\in F}I_\eta$
and all $k,l\leq\eta(\beta)$ with $k\not=l$,
$$q*\sigma\restriction\delta\forces
q^\prime\leq_{\eta(\delta)}q(\delta)\wedge
x[q^\prime*(\sigma(\delta)\frow k)]\perp
x[q^\prime*(\sigma(\delta)\frow l)].$$ Let $r\leq_{F,\eta}q$ be such
that $r\restriction\delta\forces r(\delta)=q^\prime$ and for all
$\sigma\in\prod_{\gamma\in F}I_\eta$ and all $k\leq\eta(\beta)$,
$r*\sigma\restriction\delta$ decides $x[r(\delta)*(\sigma(\delta)\frow
k)]$.

Note that $r$ is indeed $(F,\eta^\prime)$-faithful since we assumed
$q\restriction\delta$ to force that $T_{q(\delta)}$ is
$c_{\min}$-homogeneous of color~$i$.

If $\beta\not=\delta$, the argument will be similar to the one used
for Lemma \ref{limitstep}.  Let $n:=\eta(\beta)$.

For all $k\leq n$ and all $\sigma\in\prod_{\gamma\in
F}I_{\eta(\gamma)}$ let $q_{\sigma,k}$ be a name for a condition such
that
$$q*\sigma\restriction\delta\forces q_{\sigma,k}\leq
q(\delta)*\sigma(\delta)$$ and for all $l\leq n$ with $l\not= k$
$$q*\sigma\restriction\delta\forces x[q_{\sigma,k}(\delta)]\perp
x[q_{\sigma,l}(\delta)].$$

Now fix $q^\prime\leq_{F,\eta} q$ such that for all
$\sigma\in\prod_{\gamma\in F}I_{\eta(\gamma)}$ and all $k\leq n$,
$q^\prime*\sigma\restriction\delta$ decides $x[q_{\sigma,k}]$.  Note
that for all $k,l\leq n$ with $k\not= l$ we have that
$$q^\prime*\sigma\restriction\delta\forces
\parity(x[q_{\sigma,k}],x[q_{\sigma,l}])=i$$ since $[T_{q(\delta)}]$
was forced to be $c_{\min}$-homogeneous.

Choose $r$ such that $r\restriction\delta=q^\prime\restriction\delta$
and for all $\sigma\in\prod_{\gamma\in F}I_{\eta^\prime{\gamma}}$
$$r*\sigma\restriction\delta\forces
r(\delta)*\sigma(\delta)=q_{\sigma,k}$$ where $k=\sigma(\beta)(n)$
(i.e., $k$ is the last digit of $\sigma(\beta)$).

It follows from the definition of $r$ that $r\leq_{F,\eta} q$. It is
easily checked that $r$ is $(F,\eta^\prime)$-faithful.
\end{proof}

Combining the last two lemmas, we can show

\begin{lemma}\label{groundmodelsets}
Let $G$ be $\mathbb P_{\omega_2}$-generic over the ground model
$M$. Then in $M[G]$, $2^\omega$ is covered by $c_{\min}$-homogeneous
sets coded in the ground model.  In particular, in $M[G]$, $2^\omega$
is covered by $\aleph_1$ $c_{\min}$-homogeneous sets.
\end{lemma}

\begin{proof}
We work in~$M$.  Let $\dot{x}$ be a name for an element of $2^\omega$.
We show that $\dot{x}$ is forced to be a branch through a
$\parity$-homogeneous tree in~$M$.  We may assume that for some
$\alpha<\omega_2$, $\dot{x}$ is an $\mathbb P_\alpha$-name for a real
not added in a proper initial stage of the iteration $\mathbb
P_\alpha$. Clearly, $\cf(\alpha)\leq\aleph_0$.  Let $p\in\mathbb
P_\alpha$. If $\alpha$ is a limit ordinal, using Lemma
\ref{densesets}, we can decrease $p$ such that for some $i\in 2$,
$p\in E_i$.  If $\alpha$ is a successor ordinal, say
$\alpha=\delta+1$, we can use Lemma \ref{onereal} to decrease $p$ such
that for some $i\in 2$
$$p\restriction\delta\forces\mbox{``$[T_{p(\delta)}]$ is
$c_{\min}$-homogeneous of color $i$''}.$$

By induction, we define a sequence $(p_n,F_n,\eta_n)_{n\in\omega}$
such that
\begin{enumerate}
\item for all $n\in\omega$, $p_n\in\mathbb P_\alpha$, $p_n\leq p$,
$F_n\in[\alpha]^{<\aleph_0}$, $\eta_n:F_n\to\omega$, and $p_n$ is
$(F_n,\eta_n)$-faithful,
\item for all $n\in\omega$, $F_n\subseteq F_{n+1}$,
$p_{n+1}\leq_{F_n,\eta_n}p_n$, and for all $\gamma\in F_n$ we have
$\eta_n(\gamma)\leq\eta_{n+1}(\gamma)$, and
\item for all $n\in\omega$ and all $\gamma\in\supt(p_n)$ there is
$m\in\omega$ such that $\gamma\in F_m$ and $\eta_m(\gamma)\geq n$.
\end{enumerate}

This construction can be done using parts a) and b) of Lemma
\ref{limitstep} and Lemma \ref{successorstep} respectively, depending
on whether $\alpha$ is a limit ordinal or not, to extend $F_n$ or to
make $\eta_n$ bigger, together with some bookkeeping to ensure 3.  Now
$(p_n)_{n\in\omega}$ is a fusion sequence.  Let $q$ be the fusion of
this sequence.  For each $n\in\omega$ let $T_n$ be the tree generated
by $\{x[p_n*\sigma]:\sigma\in\prod_{\gamma\in
F_n}I_{\eta(\gamma)}\}$. It is easily seen that
$T_{q}=\bigcup_{n\in\omega}T_n$.

It now follows from the faithfulness of the $p_n$ that $[T_{q}]$ is
$c_{\min}$-homogeneous of color~$i$.  Clearly, $q$ forces $\dot{x}$ to
be a branch through~$T_q$.  It follows that the set of conditions in
$\mathbb P_\alpha$ forcing $\dot{x}$ to be an element of a
$c_{\min}$-homogeneous set coded in $M$ is dense in $\mathbb
P_\alpha$. Since $\mathbb P_\alpha$ is completely embedded in $\mathbb
P_{\omega_2}$, this finishes the proof of the lemma.
\end{proof}

\begin{corollary}
It is consistent with ZFC that $2^{\aleph_0}=\aleph_2$ and $2^\omega$
is covered by $\aleph_1$ $c_{\min}$-homogeneous sets, but it is not
covered by less than $2^{\aleph_0}$ $c_{\max}$-homogeneous sets.
\end{corollary}

\subsection{Why forcing with $\mathbb P_{c_{\max}}$?} One may ask
whether there is an essentially simpler way of increasing
$\hm(c_{\max})$ while keeping $\hm$ small other than iterating our
basic forcing notion $\mathbb P$. Zapletal \cite{zapletal} showed
that in certain cases there is an optimal way of increasing a
covering number of a $\sigma$-ideal.  He observed that there is
an optimal way of increasing $\hm$ in the sense that all cardinal
invariants which are not bigger than $\hm$ in ZFC are kept small
(assuming the existence of some large cardinals).  The natural
forcing to do this is the following:

\begin{defn}\label{cmindef}    % mg may 13
The $c_{\min}$-forcing $\mathbb
P_{c_{\min}}$ is the partial order consisting of all perfect
subtrees $p$ of $2^{<\omega}$ with the property that for all
$s\in p$ there are splitting nodes $t_0$ and $t_1$ of $p$ which
extend $s$ such that the length of $t_0$ is even and the length
of $t_1$ is odd.
\end{defn}

It is easy to see that the $\mathbb P_{c_{\min}}$-generic real avoids
all the $c_{\min}$-homogeneous sets in the ground model.  Therefore
iterating $\mathbb P_{c_{\min}}$ increases~$\hm$.

The natural approach for increasing $\hm(c_{\max})$ would be forcing
with an iteration of the Borel subsets of $2^\omega$ modulo the
$\sigma$-ideal generated by the $c_{\max}$-homogeneous subsets.
However, this attempt must fail.  Zapletal observed that this forcing
notion is not homogeneous, that is, the forcing notion does not stay
the same when restricted to some Borel set not covered by countably
many $c_{\max}$-homogeneous sets.  We show that in fact, this forcing
notion increases~$\hm$.

\begin{thm}
Let $X$ be any Polish space with some nontrivial continuous pair-coloring
$c:[X]^2\to 2$.  Then the Boolean algebra of Borel subsets of $X$
modulo the $\sigma$-ideal generated by the $c$-homogeneous sets is
forcing equivalent to $\mathbb P_{c_{\min}}$.
\end{thm}

The theorem easily follows from the next lemma, which is a strengthening
of Lemma \ref{canmin}.

\begin{lemma}
Assume that $A\subseteq X$ is analytic.  If $A$ is not covered by
countably many $c$-homogeneous sets, then $c_{\min}\leq c\restriction
A$, i.e., $A$ has a perfect subset on which $c$ is isomorphic to
$c_{\min}$.
\end{lemma}

\begin{proof} Since $A$ is analytic, there is a continuous map
$f:\omega^\omega\to A$ which is onto.  For $s\in\omega^{<\omega}$ let
$O_s:=\{x\in\omega^\omega:s\subseteq x\}$. For
$B\subseteq\omega^\omega$ let
\begin{multline*}B^\prime:=B\setminus\bigcup\{O_s:s\in\omega^{<\omega}\wedge
f[B\cap O_s]\mbox{ is not covered}\\\mbox{by countably many
$c$-homogeneous sets}\}.
\end{multline*}
Note that $B^\prime$ is closed if $B$ is.

Let $B_0:=\omega^\omega$, $B_{\alpha+1}:=B_\alpha^\prime$ for
$\alpha<\omega_1$ and $B_\delta:=\bigcap_{\alpha<\delta}B_\alpha$ for
limit ordinals $\delta<\omega_1$. Since there are only countably many
$O_s$, there is $\alpha<\omega_1$ such that $B_\alpha=B_{\alpha+1}$.
Let $B:=B_\alpha$.

Since $A$ is not covered by countably many $c$-homogeneous sets, $B$
is not empty.  Clearly, for every open set $O\subseteq\omega^\omega$,
$O\cap B$ is empty or $f[O\cap B]$ is not covered by countably many
$c$-homogeneous sets and therefore is not homogeneous.  It now follows
from the continuity of $f$ and $c$ that for all $s\in T(B)$ and all
$i\in 2$ there are $s_0,s_1\in T(B)$ extending $s$ such that $c$ is
constant on $f[O_{s_0}\cap B]\times f[O_{s_1}\cap B]$ with value~$i$.

This is sufficient to construct inductively a perfect binary subtree
$T$ of $T(B)$ such that $f\restriction [T]$ is 1-1 and $f[[T]]$ has
the desired properties.
\end{proof}

\section{Consistency of
$\cov(\cont(\R))<\cov(\Lip(\R))$}\label{lipcontsep} This section is
devoted to the proof of

\begin{thm} $\cov(\Cont)<\hm$ is consistent.
\end{thm}

In Definition \ref{cmindef} we have already
introduced the forcing notion $\mathbb P_{c_{\min}}$
as the right tool to increase~$\hm$.

In this section we write $\mathbb P$ for $\mathbb P_{c_{\min}}$.  As
usual, for every ordinal $\alpha$, $\mathbb P_\alpha$ denotes the
countable support iteration of $\mathbb P$ of length $\alpha$.  We
have to show

\begin{lemma}\label{covbycont} After forcing with $\mathbb P_{\omega_2}$
over a model of CH the continuous functions coded in the ground model
cover $(2^\omega)^2$ (in the extension).
\end{lemma}

How do we construct the required continuous mappings in the ground
model?  Of course, every condition $p\in\mathbb P$ is a perfect
(binary) tree and thus $[p]$ is homeomorphic to $2^\omega$.  This
homeomorphism is unique if we assume that it preserves the
lexicographic order.

Let $\alpha$ be an ordinal and $\dot{x}$ a $\mathbb P_\alpha$-name for
an element of $2^\omega$ which is not added in a proper initial stage
of the iteration.  Then for every $p\in\mathbb P_\alpha$ we construct
$q\leq p$ such that for $S:=\supt(q)$ the following  property
$(*)_{q,S,\dot{x}}$ holds:

\begin{itemize}\item[$(*)_{q,S,\dot x}$]
Let $T_q(\dot{x})$ be the tree of $q$-possibilities for $\dot{x}$
defined as in the proof of Lem\-ma \ref{onereal}.  Then in the ground
model we have a homeomorphism $h:[T_q(\dot{x})]\to(2^\omega)^S$ such
that if $G$ is $\mathbb P_\alpha$-generic with $q\in G$, then $h$ maps
$\dot{x}_G$ to a sequence $(z_\gamma)_{\gamma\in S}\in(2^\omega)^S$
such that for all $\gamma\in S$, $z_\gamma$ is the image of the
$\gamma$'th generic real under the natural homeomorphism from
$[q(\gamma)_G]$ to $2^\omega$.
\end{itemize}

So in a weak sense we can reconstruct the restriction of the sequence
of generic reals to $\supt(q)$ from $\dot{x}_G$ using a ground model
function.  We will see soon that we can really reconstruct the
sequence of generic reals below $\alpha$ from $\dot{x}_G$.

It is not difficult to see
\begin{claim}\label{decode} If $(*)_{q,S,\dot x}$
holds for some $q\in\mathbb P_\alpha$
and $S\in[\alpha]^{\leq\aleph_0}$, then also $(*)_{r,S,\dot x}$
holds  for every $r\leq q$ (with the original set $S$).
\end{claim}

Now suppose $\dot{x}$ and $\dot{y}$ are $\mathbb P_{\omega_2}$-names
for elements of $2^\omega$.  Assume that both, $\dot{x}$ and
$\dot{y}$, are forced to be new reals.  We may do so because the
constant functions take care about covering pairs $(x,y)\in
(2^\omega)^2$ where $x$ or $y$ are in the ground model.

We may also assume that there are $\alpha,\beta<\omega_2$ such that
$\dot{x}$ is in fact a $\mathbb P_\alpha$-name forced not to be added
in a proper initial stage of the iteration $\mathbb P_\alpha$ and the
same is true for $\dot{y}$ with respect to $\beta$.  Finally we may
assume $\beta\leq\alpha$.

Now let $p\in\mathbb P_\alpha$.  We find $q\in\mathbb P_\beta$ such
that $q\leq p\restriction\beta$ and $(*)_{q,\supt(q),\dot x}$ holds.
Then we can find $r\in\mathbb P_\alpha$ such that $r\leq q\frow
p\restriction [\beta,\alpha)$ and $(*)_{r,\supt(r),\dot x}$ holds.

Let $h:[T_r(\dot{x})]\to(2^\omega)^{\supt(r)}$ be the homeomorphism
(in the ground model) guaranteed by $(*)_{r,\supt(r),\dot x}$.
 Let $g:[T_r(\dot{y})]\to(2^\omega)^{\supt(q)}$ be the
homeomorphism guaranteed by $(*)_{r,\supt(q),\dot y}$,
which holds by Claim \ref{decode}.

Now let $\pi:(2^\omega)^{\supt(r)}\to(2^\omega)^{\supt(q)}$ be the
natural projection and put $f:=g^{-1}\circ\pi\circ h$.  $f$ is only
defined on a closed subset of $2^\omega$, but as in the proof of Lemma
\ref{context}, we can continuously extend it to all $2^\omega$.
Clearly $r\forces f(\dot{x})=\dot{y}$.  This finishes the proof of
Lemma \ref{covbycont} provided we know

\begin{lemma}\label{star}
Let $\alpha$ be an ordinal and $\dot{x}$ a $\mathbb P_\alpha$-name for
an element of $2^\omega$ which is not added in a proper initial stage
of the iteration.  Then for every $p\in\mathbb P_\alpha$ there is
$q\leq p$ such that $(*)_{q,\supt(q),\dot x}$ holds.
\end{lemma}

\begin{proof}  We follow closely the proof of Lemma \ref{groundmodelsets}.
We fix $\dot{x}$ throughout the following proof.

For $p\in\mathbb P$ and $n\in\omega$ let $p^n$ denote the set of those
splitting nodes of $p$ that have exactly $n$ splitting nodes among
their proper initial segments.  For $q\leq p$ we write $q\leq_n p$ if
$q^n=p^n$.  Every $\rho\in 2^n$ determines an element $s_\rho$ of
$p^n$.  Let $p*\rho:=p_{s_\rho}=\{s\in p:s\subseteq s_\rho\vee
s_\rho\subseteq s\}$.

We call a condition $p\in\mathbb P$ {\em normal} if for all splitting
nodes $s,t\in p$ such that $s\subsetneq t$ and $t$ is a minimal
splitting node above~$s$, $\dom(t)\setminus\dom(s)$ is odd, i.e.,
$\dom(s)$ and $\dom(t)$ have a different parity.

As in the $\mathbb P_{c_{\max}}$-case, if $p$ is a normal condition
and $q\leq_n p$, then there is a normal condition $r\leq_n q$.  This
is normalization at $n$ that from now on will be done automatically,
just as in the $\mathbb P_{c_{\max}}$-case.

We extend the notion of normality to conditions in $\mathbb P_\alpha$
and for $F\in[\alpha]^{<\aleph_0}$ and $\eta:F\to\omega$ we define
$\leq_{F,\eta}$ on $\mathbb P_\alpha$ as for $\mathbb P_{c_{\max}}$
(see section~\ref{forcingnotation}).
Fusion sequences are defined as for $\mathbb P_{c_{\max}}$ and it
should be clear that fusions of fusion sequences in $\mathbb P_\alpha$
are again conditions, provided the elements of the fusion sequence are
normal.

For $f$ and $\eta$ as above, $p\in\mathbb P_\alpha$, and
$\sigma\in\prod_{\gamma\in F}2^{\eta(\gamma)}$, $p*\sigma$ is defined
as in Section \ref{hmsmall}.

We also use the notion of faithfulness, but in the present context the
definition is a weaker than in Section \ref{hmsmall}.

\begin{defn} For $F$ and $\eta$ as above, $p\in\mathbb P_\alpha$ is
{\em $(F,\eta)$-faithful} iff for all $\sigma,\tau\in\prod_{\gamma\in
F}2^{\eta(\gamma)}$ with $\sigma\not=\tau$, $x[p_\sigma]\perp
x[p_\tau]$.
\end{defn}

The corresponding statement to Lemma \ref{limitstep} and Lemma
\ref{successorstep} is

\begin{claim}\label{faithfulext}
Let $F$ and $\eta$ be as before and suppose that $q\in\mathbb
P_\alpha$ is $(F,\eta)$-faithful.

a) Let $\beta\in\alpha\setminus F$ and let $F^\prime:=F\cup\{\beta\}$
and $\eta^\prime:=\eta\cup\{(\beta,0)\}$.  Then there is $r\leq_{F,\eta}q$
such that $r$ is
$(F^\prime,\eta^\prime)$-faithful.

b) Let $\beta\in F$ and let $\eta^\prime:=\eta\restriction
F\setminus\{\beta\}\cup\{(\beta,\eta(\beta)+1)\}$.  Then there is
$r\leq_{F,\eta}q$ such that $r$ is $(F,\eta^\prime)$-faithful.
\end{claim}

\begin{proof}
In contrast to the $\mathbb P_{c_{\max}}$-case, a) is not trivial
here.  This is because $\leq_0$ is not equivalent to~$\leq$.  But this
is rather a notational issue.  a) clearly follows from the proof of
b).

For b) let $\delta:=\max F$ and let $\{\sigma_0,\dots,\sigma_m\}$ be
an enumeration of $\prod_{\gamma\in F}2^{\eta(\gamma)}$.  We define a
$\leq_{F,\eta}$-decreasing sequence $(q_j)_{j\leq m}$ along with names
$q_{\sigma,0}$ and $q_{\sigma,1}$, $\sigma\in\prod_{\gamma\in
F}2^{\eta(\gamma)}$, for conditions.

Let $j\in\{1,\dots m\}$ and assume that $q_{j-1}$ has been constructed
already.  Since $\dot{x}$ is not added in a proper initial stage of
the iteration, there are $q_{\sigma_j,0}$ and $q_{\sigma_j,1}$ such
that for all $i\in 2$
$$q_{j-1}*\sigma_j\restriction\delta\forces q_{\sigma_j,i}\leq
(q(\delta)*(\sigma_j(\delta)\frow i))\frow q\restriction
(\delta,\alpha)$$ and
$$q_{j-1}*\sigma_j\restriction\delta\forces x[q_{\sigma_j,0}]\perp
x[q_{\sigma_j,1}].$$ Let $q_j\leq_{F,\eta} q_{j-1}$ be such that
$q_j*\sigma\restriction\delta$ decides $x[q_{\sigma_j,0}]$ and
$x[q_{\sigma_j,1}]$.  This finishes the inductive construction of the
$q_j$.

Now let $r\leq_{F,\eta} q_m$ be such that
$r\restriction\delta=q_m\restriction\delta$ and for all
$\sigma\in\prod_{\gamma\in F}2^{\eta(\gamma)}$ and all coordinatewise
extensions $\tau\in\prod_{\gamma\in F}2^{\eta^\prime(\gamma)}$ of
$\sigma$,
$$r*\tau\restriction\delta\forces
r*\tau\restriction[\delta,\alpha)=q_{\sigma,\tau(\eta(\beta))}.$$

It is easy to check that $r$ works for the claim.
\end{proof}

To conclude the proof of Lemma \ref{star}, let $p\in\mathbb P_\alpha$.
Using some bookkeeping and parts a) and b) of Claim \ref{faithfulext}
we construct a sequence $(p_n)_{n\in\omega}$ and a sequence
$(F_n,\eta_n)_{n\in\omega}$ witnessing that $(p_n)_{n\in\omega}$ is a
fusion sequence such that $p=p_0$ and for all $n\in\omega$, $p_n$ is
$(F_n,\eta_n)$-faithful.

Let $q$ be the fusion of the sequence $(p_n)_{n\in\omega}$.  We have
to check that $(*)_{q,\supt(q),\dot x}$
 holds.

Let $a\in[T_q(\dot{x})]$ and $n\in\omega$. Now $q\leq_{F_n,\eta_n}p_n$
and $p_n$ is $(F_n,\eta_n)$-faithful.  It follows that there is
exactly one $\sigma_{a,n}\in\prod_{\gamma\in F_n}2^{\eta_n(\gamma)}$
such that $x[q_\sigma]\subseteq a$.

Let
$h(a):=(\bigcup_{n\in\omega}\sigma_{a,n}(\gamma))_{\gamma\in\supt(q)}$.
Since for all $\gamma\in\supt(q)$ and all $m\in\omega$ there is some
$n\in\omega$ such that $\gamma\in F_n$ and $\eta_n(\gamma)\geq m$,
$h(a)\in(2^\omega)^{\supt(q)}$.  It is easily checked that
$h:[T_q(\dot{x})]\to(2^\omega)^{\supt(q)}$ is a homeomorphism
witnessing $(*)_{q,\supt(q),\dot x}$.
\end{proof}

\section{Concluding remarks and open problems}
The numbers  $\hm(c_{\min})$,
$\hm(c_{\max})$, $\cov(\Cont(\R))$ and $\cov(\Lip(\R))$ are examples of
covering numbers of meager ideals. Although the hope expressed by
Blass in \cite{blass} to find a classification of all ``simple''
cardinal invatiants of the continuum was shattered by the construction
in \cite{448} of uncountably many different covering numbers of simply
defined meager ideals, there is still hope to find the ``largest''
nontrivial covering number of a meager ideal. By ``nontrivial'' it is
meant that the number can consistently be smaller than the continuum.

At the moment the leading candidate for such a number is
$\hm(c_{\max})$.  The numbers $\cov(\Cont(\R))$ and $\cov(\Lip(\R))$
are also very large, and perhaps larger nontrivial covering numbers of
meager ideals can be found by considering covering by functions with a
stronger regularity condition than Lipschitz. It would be natural to
compare $\hm(c_{\max})$ to covering by smooth (total)
functions. At the moment nontriviality is open even for
differentiable functions. At any rate, the ideal generated by real
\emph{analytic} functions is certainly too small:  Every analytic 
function is either constant or the graph of the function intersects every 
horizontal line only in countably many points.  This easily implies that 
less than $2^{\aleph_0}$ analytic functions cannot cover $\mathbb R^2$.

The meager ideals which historically led to the study of homogeneity
numbers are the \emph{convexity ideals} of closed subsets of
$\R^2$. If a closed subset of a Euclidean space is not covered by
countably many convex subsets (namely, its convex subsets generate a
proper $\sigma$-ideal), it has a closed subset on which the convex
subsets of the whole set generate a meager ideal (see \cite{kojman} or
\cite{convexitynumbers}). For some closed subsets of the plane, this
meager ideal coincides with the homogeneity ideal of some continuous
pair coloring \cite{GKKS}.

Saharon Shelah remarked recently to the authors that he came close to
discovering the properties of $\hm$ in his investigations of monadic
theory of order \cite{shelahpersonal}. In an attempt to remove GCH
from the proof in the last section of \cite{monadic} Shelah found a
proof from the assumption $\hm=2^{\aleph_0}$.  He was able to prove
that $\hm=2^{\aleph_0}$ if the continuum is a limit cardinal, but did
not prove more about $\hm$ and eventually found a way to eliminate GCH
which did not involve $\hm$, which was consequently published in
\cite{GurevichShelah}.

 It is not clear why homogeneity numbers of continuous pair-colorings
on Polish spaces were not studied earlier. We can only speculate about
that. In the very short time since their study was begun, these
numbers were related to quite a few subjects. Apart from the relation
to planar convex geometry and to finite random graphs, which were
mentioned above, there are relations to large cardinals, determinacy
and pcf theory. Quite recently, Shelah and Zapletal \cite{duality}
defined $n$-dimensional generalizations of $\hm(c_{\min})$ and
integrated forcing, pcf theory and determinacy theory to prove several
duality theorem for those numbers.

We do not know at the moment if $\aleph_1<\hm<2^{\aleph_0}$ is
consistent or not. We do not know if there is a closed planar set
whose convexity number is equal to $\hm(c_{\max})$. We also find
the following intriguing:

\begin{problem}
Are the equalities $\hm(c)=\hm(c_{\min})$
and $\hm(c)=\hm(c_{\max})$ which hold in $V^{P_{c_{\max}}}$ absolute for
a reduced coloring~$c$ on a Polish space $X$?  In other words, does a
reduced coloring $c$ that satisfies $\hm(c)=\hm(c_{\max})$ in some
model of set theory which separates $\hm(c_{\min})$ and $\hm(c_{\max})$
satisfy this in every model that separates
$\hm(c_{\min})$ and $\hm(c_{\max})$?
\end{problem}

\end{document}